\documentclass[10pt]{article} % For LaTeX2e
\usepackage[accepted]{tmlr}
\usepackage{amsmath,amsfonts,bm}

\def\eqref#1{equation~\ref{#1}}
\def\1{\bm{1}}

\def\vtheta{{\bm{\theta}}}

\DeclareMathAlphabet{\mathsfit}{\encodingdefault}{\sfdefault}{m}{sl}
\SetMathAlphabet{\mathsfit}{bold}{\encodingdefault}{\sfdefault}{bx}{n}

\newcommand{\KL}{D_{\mathrm{KL}}}

\DeclareMathOperator*{\argmin}{arg\,min}

\usepackage{microtype}
\usepackage{graphicx}
\usepackage{xcolor}
\usepackage{subcaption}
\usepackage{enumitem}
\usepackage{booktabs} % for professional tables
\usepackage[dvipsnames]{xcolor}
\usepackage{comment}
\usepackage{makecell}

\usepackage{hyperref}
\usepackage{url}
\usepackage{float}
\usepackage{tcolorbox}
\newtcolorbox{roundedbox}[1][]{%
  colback=white, % 无背景色
  colframe=black, % 边框颜色为黑色
  arc=3pt, % 圆角大小
  boxrule=0.8pt, % 边框粗细
  width=\linewidth, % 宽度
  boxsep=0.0pt, % 内部间距
}

\usepackage[normalem]{ulem}

\title{Convergence Analysis of Wasserstein Proximal Algorithm beyond Geodesic Convexity}

\author{\name Shuailong Zhu \email shzhu@math.ubc.ca \\
      \addr Department of Mathematics \\
      University of British Columbia
      \AND
      \name Xiaohui Chen \email xiaohuic@usc.edu \\
      \addr Department of Mathematics \\
      Thomas Lord Department of Computer Science \\
      University of Southern California \\ \\
      Computing and Mathematical Sciences \email xiaohuic@caltech.edu \\
      California Institute of Technology}

\usepackage{amsmath}
\usepackage{amssymb}
\usepackage{mathtools}
\usepackage{amsthm}
\usepackage[capitalize,noabbrev]{cleveref}
\theoremstyle{plain}
\newtheorem{theorem}{Theorem}[section]

\newtheorem{lemma}[theorem]{Lemma}
\newtheorem{corollary}[theorem]{Corollary}
\newtheorem{example}[theorem]{{Example}}
\theoremstyle{definition}
\newtheorem{definition}[theorem]{Definition}
\newtheorem{assumption}{Assumption}
\theoremstyle{remark}
\newtheorem{remark}[theorem]{Remark}

\def\month{09}  % Insert correct month for camera-ready version
\def\year{2026} % Insert correct year for camera-ready version
\def\openreview{\url{https://openreview.net/forum?id=GU9JPTcBsM}} % Insert correct link to OpenReview for camera-ready version

\begin{document}

\maketitle
\begin{abstract}
  The proximal algorithm is a powerful tool to minimize nonlinear and nonsmooth functionals in a general metric space. Motivated by the recent progress in studying the training dynamics of the noisy gradient descent algorithm on two-layer neural networks in the mean-field regime, we provide in this paper a simple and self-contained analysis for the convergence of the general-purpose Wasserstein proximal algorithm without assuming geodesic convexity on the objective functional. Under a natural Wasserstein analog of the Euclidean Polyak-{\L}ojasiewicz inequality, we show that the {proximal algorithm} achieves an \emph{unbiased} and \emph{dimension-free} linear convergence rate {at the population level}. Our convergence rate improves upon existing rates of the proximal algorithm for solving Wasserstein gradient flows when specialized to strong geodesic convex functionals. We also extend our analysis to the inexact proximal algorithm for geodesically semiconvex objectives. In our numerical experiments, proximal training demonstrates a faster convergence rate than the noisy gradient descent algorithm on two-layer mean-field neural networks\footnote{Code is available at \href{https://github.com/One-punch24/Convergence-Analysis-of-WPA}{https://github.com/One-punch24/Convergence-Analysis-of-WPA}.}.
\end{abstract}

\section{Introduction}
Minimizing a cost functional over the space of probability distributions has recently drawn widespread statistical and machine learning applications such as variational inference~\citep{lambert2022variational,ghosh2022,yao2023meanfield}, sampling~\citep{wibisono18a,ULA2019,chewi2024analysis}, and generative modeling~\citep{xu2024normalizing,cheng2024}, among many others. In this work, we consider the following general optimization problem:
\begin{equation}
\label{Eqn_objective}
    \min_{\rho\in\mathcal{P}_2({\Theta})}F(\rho),
\end{equation} %and $\rho^*$ is defined as a solution to (\ref{Eqn_objective}).
where $F$ is a real-valued target functional defined on the space of probability distributions $\mathcal{P}_2({\Theta})$ with finite second moments on a domain $\Theta\subset \mathbb{R}^d$. Our motivation for studying this problem stems from analyzing training dynamics of the Gaussian noisy gradient descent algorithm on infinitely wide neural networks, which can be viewed as a forward time-discretization of the mean-field Langevin dynamics (MFLD)~\citep{mei2019meanfield, hu2021mean}. Given the connection between sampling and optimization, the continuous-time MFLD is an important example of the Wasserstein gradient flow corresponding to minimizing an entropy-regularized total objective function of large interacting particle systems (cf. Section~\ref{subsec:WGF_LD}).

On the other hand, the Wasserstein gradient flow is conventionally constructed by the following proximal algorithm,
\begin{equation}
\label{Eqn_jko}
    \rho_{n+1} \in prox_{F,\xi}(\rho_n) := \argmin_{\tilde{\rho}\in\mathcal{P}_2({\Theta})} \Big\{ F(\tilde{\rho})+\dfrac{1}{2\xi}\mathcal{W}_2^2(\tilde{\rho},\rho_n) \Big\},
\end{equation}
where $\xi > 0$ is the time-discretization step size. Here, we do not require the proximal operator to have a unique minimizer---our main convergence results in Section~\ref{Sec_main} for the Wasserstein proximal algorithm hold for any minimizer in $prox_{F,\xi}(\rho_n)$. The Wasserstein proximal algorithm (\ref{Eqn_jko}) is an iterative \emph{backward time-discretization} procedure for solving (\ref{Eqn_objective}), and it is also known as the JKO scheme introduced in the seminal work~\citep{jko}. In contrast to various forward-discretization methods such as gradient descent over $\mathcal{P}_2({\Theta})$ and the Langevin sampling algorithms~\citep{durmus2019analysis,ULA2019,chewi2024analysis}, proximal algorithms are often {\it unbiased} in the sense that their convergence guarantees do not depend on the dimension-dependent discretization error with positive step size and they are often more stable than the forward gradient descent algorithms without strong smoothness condition~\citep{yao2023meanfield,yao2024wasserstein,salim2020wasserstein,cheng2024}. While the proximal algorithms for geodesically convex functionals are well studied {in the literature \citep{naldi2021weak,di2025inexact}}, it remains an open question 
\emph{whether they can maintain similar unbiased and linear convergence guarantees in discrete time beyond the geodesic convexity}. Current work fills this important gap by establishing linear convergence results without assuming geodesic convexity on the objective functional $F$.

Our general quantitative convergence rate for the Wasserstein proximal algorithm offers an alternative training scheme to the noisy gradient descent for two-layer neural networks in the mean-field regime. Specifically, a two-layer neural network is parameterized as
\begin{equation}
\label{eq:twolayer_nn}
    f(x ; \vtheta) := \dfrac{1}{m}\sum_{j=1}^m \varphi(\theta_j;x)=\int_{\Theta} \varphi(\theta;x)d\rho^m(\theta),
\end{equation}
where $\vtheta =(\theta_1,\theta_2,...,\theta_m)\in \mathbb{R}^{d\times m}$ and $\rho^m=\frac{1}{m}\sum_{j=1}^m \delta_{\theta_j}$ is the empirical distribution of the hidden neuron parameters. The perceptron $\varphi(\theta_j ; x)$ in (\ref{eq:twolayer_nn}) can take the form $\varphi(\theta_j ; x) = \sigma(\theta_j^\top x)$ where $\sigma$ is some nonlinear activation function. Given a training dataset $(x_i,y_i)_{i=1}^N \sim p(x,y)$ and a convex loss function $l(\cdot)$ (such as the  squared loss and logistic loss), the $L^2$-regularized training risk is defined as 
\begin{align}
\label{twonn}
   \begin{gathered}
    % &= {1 \over N} \sum_{i=1}^N l( f(x_i ; \vtheta),y_i ) +\dfrac{\lambda}{m}\sum_{j=1}^m \left\Vert \theta_j\right\Vert^2\\
       R({\rho^m})= {1 \over N} \sum_{i=1}^N l \left( \int_{\Theta} \varphi(\theta ; x_i)d\rho^m(\theta),y_i \right)
       +\lambda \int_{\Theta} \left\Vert\theta\right\Vert^2 d\rho^m(\theta),
   \end{gathered}
\end{align}
 where $\lambda>0$ is the coefficient of $L_2$-regularization. The Gaussian noisy gradient descent algorithm on the $L^2$-regularized training risk can be written as the following stochastic recursion
\begin{equation}
\label{noise_gd_eqn}
    \theta^{n+1}_j = \theta^{n}_j - \xi\nabla {\delta R \over \delta \rho} (\frac{1}{m}\sum_{\ell=1}^m \delta_{\theta_\ell^n}) (\theta^{n}_j) + \sqrt{2\xi\tau }z_{nj},
\end{equation}
where $z_{nj}$ are i.i.d. $\mathcal{N}(0,I)$ and $\tau>0$ represents the Gaussian noise variance. In (\ref{noise_gd_eqn}), ${\delta R \over \delta \rho}(\rho)$ is the first variation of $R$ at $\rho$ (cf. Definition~\ref{def:first_variation}). The limiting dynamics of (\ref{noise_gd_eqn}) under $m\rightarrow \infty$ and $\xi\rightarrow 0$ is called the continuous-time MFLD \citep{hu2021mean}. Under a uniform log-Sobolev inequality (LSI) assumption (cf. Definition~\ref{def:unifLSI}), linear convergence of MFLD to the optimal value of the total objective ($L^2$-training risk plus an entropy term) is established in \cite{nitanda2022convex,chizat2022meanfield}, and the noisy gradient descent algorithm is subject to a \emph{dimension-dependent} time-discretization error~\citep{nitanda2022convex}, which may slow down the convergence.

To remove the time-discretization error, we may instead train the neural network with the Wasserstein proximal algorithm (\ref{Eqn_jko}). Since such neural network architecture satisfies the uniform LSI which in turn implies a Wasserstein Polyak-{\L}ojasiewicz (PL) inequality (cf.  Definition \ref{Def_PL}), the exact algorithm can theoretically achieve an unbiased linear rate of convergence to a global minimum of the total objective {at the population level}.

\subsection{Contributions}
  In this work, we give a simple and self-contained convergence rate analysis of the Wasserstein proximal algorithm (\ref{Eqn_jko}) for minimizing the objective function satisfying a (Wasserstein) PL-type inequality (\ref{Eqn_wass_PL}) without resorting to any geodesic convexity assumption. %In the MFLD setting, our PL inequality relaxes the standard uniform log-Sobolev inequality. 
  Below we summarize our main contributions.
  \begin{itemize}[noitemsep,nolistsep,leftmargin=*]
      \item  To the best of our knowledge, current work is among the first works to obtain a {theoretically} \emph{unbiased} and \emph{dimension-free} linear convergence rate of the {exact} Wasserstein proximal algorithm for optimizing a functional under \textbf{merely a PL-type inequality}. Our analysis applied to $\mu$-convex ($\mu>0$) objective functional along geodesics   yields a faster linear convergence rate than the existing literature \citep{yao2023meanfield,cheng2024}.

      \item The linear convergence guarantee provides a new training scheme for two-layer wide neural networks in the mean-field regime. Our numerical experiments show a faster training phase (up to particle discretization error) than the (forward) noisy gradient descent method.
      
      \item We also analyze the inexact proximal algorithm for geodesically semiconvex objectives under the PL-type inequality.
  \end{itemize}
{The paper \citep{lascu2024linear} also proves the
exponential convergence under flat convexity and a uniform LSI for several JKO-based schemes, including a fully implicit proximal-point scheme as well as prox-linear and proximal-gradient schemes. The main difference between this and our work is that our general convergence analysis is formulated directly under the Wasserstein PL inequality, without requiring either flat convexity or geodesic convexity while \cite{lascu2024linear} focuses on MFLD setting. In particular, Example \ref{Cor_specialexample} in our paper provides an explicit functional that satisfies the Wasserstein PL inequality while being neither flatly convex nor geodesically convex.  These two works were conducted independently.}

\subsection{Literature review}
\label{Sec_review}
Recently, various time-discretization methods have been proposed for minimizing a functional over a single distribution. Different from the proximal algorithm, some explicit forward schemes that can be seen as gradient descent in Wasserstein space are proposed \citep{chewi2020gradient,liu2016stein}. For example, \cite{chewi2020gradient} studies a gradient descent algorithm for solving the barycenter problem on the Bures-Wasserstein manifold of Gaussian distributions.  The Langevin algorithm, as another forward discretization of Wasserstein gradient flows via its stochastic differential equation (SDE) recursion, is widely used in the sampling literature. %\footnote {Time-discretization of SDEs, e.g. Langevin dynamics (\ref{Eqn_langevin}), MFLD (\ref{Eqn_mfld}), are all called Langevin algorithm in our paper}. 
Numerous works \citep{durmus2019analysis,ULA2019,chewi2024analysis} have been devoted to the analysis of the Langevin algorithm under different settings and its variants~\citep{zhang2023improved,yuansi2022}. However, Langevin algorithms are naturally biased for a positive step size. \cite{salim2020wasserstein} introduced a hybridized forward-backward discretization, namely the Wasserstein proximal gradient descent, and proved convergence guarantees for geodesically convex objectives, akin to the proximal gradient descent algorithm in Euclidean spaces.  {\citet{luu2024non} proposed a semi Forward-Backward scheme to optimize a subclass of non-geodesically convex objectives with a difference-of-convex structure.}

\textbf{Existing rate analysis for proximal algorithm.} Though convergence rate analysis for Langevin algorithms under strong convexity is well-developed, it is not until recently that the convergence rate of the proximal algorithm on geodesically convex objectives is obtained. One advantage of the proximal algorithm is that it ensures a dimension-independent convergence guarantee directly for any starting distribution. \cite{yao2023meanfield,cheng2024} proved an unbiased linear convergence result for the $\mu$-strongly convex objective. The condition is relaxed to geodesic convexity and quadratic growth of functional in \citep{yao2024wasserstein}. {\cite{di2025inexact} provided a thorough analysis of inexact proximal algorithm under convexity along generalized geodesics.} However, convergence analysis of proximal algorithm for non-geodesically convex objective functionals remains largely missing. 

\textbf{Convergence rate of different time-discretizations under PL-type inequality.} \cite{ULA2019} obtained a biased linear convergence result for Langevin dynamics under the log-Sobolev inequality (LSI) and smoothness condition. \cite{nitanda2022convex} extended this result to MFLD with similar techniques. The proximal Langevin algorithm proposed by {\cite{wibisono2019proximallangevinalgorithmrapid}, attains a biased linear convergence rate under the LSI, while an extra smoothness condition of the second derivative of the sampling function is required. Proximal sampling algorithm~\citep{cyx2022}, assuming access to samples from an oracle distribution, achieves an unbiased linear convergence for sampling from Langevin dynamics under the LSI, while the analysis requires geodesic semiconvexity (cf. Definition~\ref{def:mu_conxity}). \cite{fanjiao23a,liang2024proximaloraclesoptimizationsampling} improved the results, however, their focus is still on sampling on a fixed function and cannot be applied to MFLD. 

To highlight the distinction between our contributions and existing results from the literature, we make the following comparison in Table \ref{tb_1} between explicit convergence guarantees of the Wasserstein proximal algorithm and Langevin algorithms for optimizing the KL divergence (\ref{eqn:KL_div}). Similar comparison on the convergence rates can be made between our result and the forward time-discretization of MFLD under further assumptions~\citep{nitanda2022convex,chizat2022meanfield}. In Table~\ref{tb_1}, we would like to make a comment on the proximal convergence rates (i.e., Rows 2 and 3). Note that $\mu$-strong convexity of a general functional $F$ does not directly imply $\mu$-PL inequality, and the rate we obtained for the proximal algorithm in Row 3 of Table~\ref{tb_1} does not require the verification of Assumption~\ref{assum_2} (cf. Theorem~\ref{Cor_sharp}). Thus, despite the same convergence rate of our results, the analysis of Row 3 in our paper follows Theorem \ref{Cor_sharp} and does not directly follow from the analysis of Row 2 (Corollary~\ref{cor_ld}).

\begin{table*}[ht]
\caption{Comparison between Langevin (forward discretization) and Wasserstein proximal algorithms for the KL divergence functional with a target distribution $\nu \propto e^{-f}$.}
\label{tb_1}
\centering
\begin{tabular}{c c c c}
\hline Algorithm &Assumptions &Step size & Convergence guarantee at $n$-th iteration\\
\hline
% \centering\makecell{Langevin\\ Algorithm}&  \centering{\makecell{$\mu$-strongly convex  \\ $L$-smoothness}} & $0<\xi<\frac{} $ & Column 4 \\ \cline{1-4} 

 \makecell{Langevin} & \centering{\makecell{$\nu$ is $\mu$-LSI  \\ $f$ is $L$-smooth (on $\Theta$)}} & $0<\xi<\dfrac{\mu}{4L^2}$ & \makecell{$\underset{\displaystyle\text{\cite{ULA2019}}}{e^{-n\mu\xi}\KL(\rho_0 \| \nu)+\dfrac{8\xi dL^2}{\mu}}$} \\ \cline{1-4} 
 \makecell{Proximal}& \centering{\makecell{$\nu$ is $\mu$-LSI  \\ $f$ is semiconvex (on $\Theta$)}} & $\xi>0$ & \makecell{$\underset{\displaystyle\text{{Ours, Corollary~\ref{cor_ld}}}}{\dfrac{1}{(1+\mu\xi)^{2n}}\KL(\rho_0 \| \nu)}$} \\ \hline
 \rule{0pt}{4.5ex}\makecell{Proximal}& \centering{\makecell{$f$ is $\mu$-strongly convex \\(on $\Theta$) \\ }}  & $\xi>0$ &\makecell{$\underset{\displaystyle\text{\cite{yao2023meanfield}}}{\dfrac{1}{(1+\mu\xi)^{n}}\KL(\rho_0 \| \nu)}$ $\rightarrow\underset{\displaystyle\text{{Ours, Theorem~\ref{Cor_sharp}}}}{\dfrac{1}{(1+\mu\xi)^{2n}}\KL(\rho_0 \| \nu)}$}  \\ \hline
\end{tabular}
\end{table*}

The remainder of this paper is organized as follows. %In Section \ref{Sec_review}, we have a basic literature review for different discretization schemes for Wasserstein gradient flows, and compare their convergence rate with the proximal algorithm under different conditions. 
In Section \ref{Sec_prelim}, we provide some background knowledge for the connection between Wasserstein gradient flows and associated Langevin dynamics. In Section \ref{Sec_main}, we present our main convergence results. In Section \ref{Sec_apply}, we discuss how to apply the proximal algorithm for MFLD of a two-layer neural network and provide numerical experiments exploring the behavior of the proximal algorithm.

%\subsection{Notations}
{\bf Notations.}
{We assume $\Theta = \mathbb{R}^d$ (by default) throughout the paper unless explicitly indicating that it is a compact subset of $\mathbb{R}^d$.} Let $\mathcal{P}_2(\Theta)$ be the collection of all probability measures with finite second moment, and $\mathcal{P}_2^a(\Theta)\subset \mathcal{P}_2(\Theta)$ be the subset of absolutely continuous measures {with respect to Lebesgue measure on $\mathbb{R}^d$}. For a measurable map $T:\Theta\rightarrow \Theta$, let $T_{\#}: \mathcal{P}_2(\Theta)\rightarrow \mathcal{P}_2(\Theta)$ be the corresponding pushforward operator. For probability measures $\mu$ and $\nu$, we shall use $T_\mu^\nu$ to denote the optimal transport (OT) map from $\mu$ to $\nu$ and $\mathbf{id}$ to denote the identity map. We use $\mathcal{W}_2(\cdot, \cdot)$ to denote the 2-Wasserstein distance. We denote $\partial F(\rho)$ to be the (Fr\'echet) subdifferential at $\rho\in\mathcal{P}_2^a(\Theta)$ whenever it exists, $\mathcal D(F):=\{\rho \in \mathcal{P}_2(\Theta) \mid F(\rho)<\infty\}$ to be the domain of $F$ that has finite functional value, and $D(|\partial F|)$ to be the domain of $F$ that has finite metric slope, see Lemma 10.1.5 of \cite{ambrosio2005gradient}. {Throughout this paper, ``$\subset$'' means ``is a subset of,'' not necessarily ``is a proper subset of.''} We refer to Appendix \ref{Appendix_background} for more notions and definitions.

\section{Preliminaries}
\label{Sec_prelim}

In this section, we review the connection between Wasserstein gradient flows and the associated Langevin dynamics. %\textcolor{red}{We introduce a strong lemma from \citep{Ambrosio_2013}, which is essential in our proof of main theorem (Theorem \ref{main_thm}).} We refer to Appendix \ref{Appendix_background} for some basic notions and differential structures in optimal transport.

\subsection{Wasserstein gradient flows and continuous-time Langevin dynamics}
\label{subsec:WGF_LD}

Gradient flows in the Wasserstein space of probability distributions provide a powerful means to understand and develop practical algorithms for solving diffusion-type equations~\citep{ambrosio2005gradient}. For a smooth Wasserstein gradient flow, noisy gradient descent algorithms over relative entropy functionals are often used for space-time discretization via the stochastic differential equation (SDE). Below we illustrate two main Wasserstein gradient flow examples involving the linear and nonlinear Fokker-Planck equations, which model the diffusion behavior of probability distributions.

{\bf Langevin dynamics via the Fokker-Planck equation.} The Langevin dynamics for the target distribution $\nu\propto e^{-f}$ is defined as an SDE,
\begin{equation}
\label{Eqn_langevin}
    d\theta_t=-\nabla f(\theta_t)dt+\sqrt{2}dW_t,
\end{equation}
where $W_t$ is the standard Brownian motion in $\Theta$ with zero initialization. It is well-known that, see e.g., Chapter 8 of~\cite{santambrogio2015optimal}, if the process $(\theta_t)$ evolves according to the Langevin dynamics in (\ref{Eqn_langevin}), then their marginal probability density distributions $\rho_t(\theta)$ satisfy the Fokker-Planck equation
\[
\partial_t \rho_t - \Delta \rho_t - \nabla \cdot (\rho_t \nabla f) = 0,
\]
which is the Wasserstein gradient flow for minimizing the KL divergence (i.e., the relative entropy)
\begin{equation}
    \label{eqn:KL_div}
    \KL(\rho \| \nu) %&:=& \int_{\Theta} \rho(\theta) \log \dfrac{\rho(\theta)}{\nu(\theta)}d\theta \\
    = \int_{\Theta} f(\theta) \rho(\theta) d \theta + \int_{\Theta} \rho(\theta) \log \rho(\theta) d\theta.
\end{equation}
% Obviously, the KL objective functional is a special case of (\ref{Eqn_objective}) as an entropy regularized linear potential functional in $\rho$. 
If $\nu$ satisfies a log-Sobolev inequality (LSI) with constant $\mu>0$, i.e., if for all $\rho\in \mathcal{P}_2^a(\Theta)$, we have
\begin{equation}
\label{eqn:LSI}
    \KL(\rho \| \nu)\leq \dfrac{1}{2\mu}J(\rho \| \nu),
\end{equation}
where $J(\rho \| \nu)=\int_{\Theta} \left\Vert \nabla \log\frac{\rho}{\nu}\right\Vert^2 d\rho$ is the relative Fisher information, then the continuous-time Langevin $\rho_t$ converges to $\nu$ exponentially fast~\citep{bakry2013analysis}.

%As $J_v(\rho)=\int \Vert \nabla \frac{\delta KL_v}{\delta \rho}(\rho)\Vert^2d\rho$, Langevin dynamics under LSI is an example that satisfies PL-type inequality (\ref{Eqn_wass_PL}), for which linear convergence of WGF can be shown.

{\bf Mean-field Langevin dynamics (MFLD) via the McKean-Vlasov equation.} In an interacting $m$-particle system, the potential energy contains a nonlinear interaction term in addition to $\int_{\Theta} f d \rho$ in (\ref{eqn:KL_div}). More generally, in the mean-field limit as $m \to \infty$, the nonlinear Langevin dynamics can be described as
\begin{equation}
\label{Eqn_nonlinearlangevin}
    d\theta_t = -\nabla {\delta R \over \delta \rho}(\rho_t)(\theta_t)dt+\sqrt{2 \tau}dW_t,
\end{equation}
where $R : \mathcal{P}_2({\Theta}) \to \mathbb{R}$ is a cost functional such as the $L^2$-regularized training risk of 
 mean-field neural networks in (\ref{twonn}) and $\tau > 0$ is a temperature parameter. For a convex loss $l$, the risk $R$ in (\ref{twonn}) has linear convexity. Process evolving according to (\ref{Eqn_nonlinearlangevin}) solves the following McKean-Vlasov equation~\citep{pmlr-v178-yao22a},
\begin{equation}
\label{eqn:McKeanVlasov}
    \partial_t \rho_t - \tau \Delta \rho_t - \nabla \cdot (\rho_t \nabla {\delta R \over \delta \rho}(\rho_t)) = 0,
\end{equation}
which is the Wasserstein gradient flow of the entropy-regularized total objective, 
\begin{equation}
\label{eqn:MFLD_obj}
    F_\tau(\rho) = R(\rho) + \tau \int_{\Theta} \rho \log \rho.
\end{equation}
Similarly, as in the linear Langevin case if the proximal Gibbs distribution of $\rho$ 
satisfies a uniform LSI (cf. Definition~\ref{def:unifLSI}), %\emph{and} the cost functional $R$ is convex in the linear structure, 
then MFLD converges to the optimal value exponentially fast in continuous time~\citep{chizat2022meanfield,nitanda2022convex} and, in the case of infinite-width neural networks in mean-field regime, it is subject to a dimension-dependent time-discretization error~\citep{nitanda2022convex}.

\section{Convergence rate analysis}
\label{Sec_main}
In this section, we first introduce a natural PL inequality in the Wasserstein space and then provide the convergence rate analysis for the Wasserstein proximal algorithm under such a weak assumption. Then, we shall extend our analysis to the inexact proximal algorithm setting. Throughout this section, we make the following regularity assumption.

\begin{roundedbox}

\begin{assumption}[Regularity assumption] The functional $F:\mathcal{P}_2(\Theta)\rightarrow (-\infty,+\infty]$ satisfies \\
(1) Proximal algorithm (\ref{Eqn_jko}) admits a minimizer for any $\rho\in D(F)$ and $\xi>0$;\\
(2)  $\mathcal{D}(F)\subset \mathcal{P}_2^a(\Theta)$. 
\label{assum_1}
\end{assumption} 
\end{roundedbox}
Assumption \ref{assum_1} ensures that the proximal algorithm admits a minimizer in $\mathcal{P}_2^a(\Theta)$, {which allows us to use Brenier's theorem and the Monge-map formulation in the proof.} Weak lower semicontinuity {and lower boundedness} of a proper functional $F$ is a sufficient condition for the existence of a minimizer (cf. Lemma \ref{Lemma_exist}), and some examples of weakly lower semicontinuous functions are provided in Remark \ref{rmk_weakly_conti}.  We refer the reader to Remark \ref{rmk_admit_minimizer} for other conditions that guarantee (1) in Assumption \ref{assum_1}.

%\subsection{A strong technical Lemma}
\begin{definition}[Hopf-Lax formula]
\label{Def_hopf}
    Let $\xi > 0$ and $\rho_\xi \in prox_{F,\xi}(\rho)$. The Hopf-Lax formula $u(\rho,\xi)$ of a functional $F : \mathcal{P}_2(\Theta) \to (-\infty,+\infty]$ is defined as
\begin{equation}
\label{eqn:Hopt-Lax}
\begin{gathered}
    u(\rho,\xi) := F(\rho_{\xi})+\dfrac{1}{2\xi}\mathcal{W}_2^2(\rho_\xi,\rho). %\\
   % \qquad \qquad = \min_{\tilde{\rho}\in\mathcal{P}_2(\Theta)} \Big\{ F(\tilde{\rho})+\dfrac{1}{2\xi}\mathcal{W}_2^2(\tilde{\rho},\rho) \Big\}.
\end{gathered}
\end{equation}
\end{definition}

{The Hopf-Lax formulation in~(\ref{eqn:Hopt-Lax}) is also known as the Moreau-Yoshida approximation \citep{ambrosio2005gradient}.} Below, we present a key connection between the time-derivative of the Hopf-Lax semigroup and the squared Wasserstein distance between the proximal and the initial point.

\begin{lemma}\label{Lemma_important}
Let $\rho\in \mathcal D(F)$. Under Assumption \ref{assum_1},
\begin{align}
    \partial_{\xi}u(\rho,\xi)&=-\frac{1}{2\xi^2}\mathcal{W}_2^2(\rho_\xi,\rho)
    \label{Eqn_important_lemma}
\end{align}
holds for $\xi\in(0,+\infty)$ with at most countable exceptions, where $\rho_\xi$ can be any proximal operator minimizer in~(\ref{eqn:Hopt-Lax}). 
\end{lemma}

\begin{remark}[Computation of the proximal operator]
    To compute $\rho_{n+1}$ in~(\ref{Eqn_jko}), we can reformulate the proximal algorithm into an optimization problem in functional space. Finding $\rho_{n+1}$ is equivalent to finding an optimal transport (OT) map $T$ such that $T_{\#}\rho_{n}$ minimizes (\ref{Eqn_jko}),
\begin{equation}
\label{Eqn_prox_op_compute}
T_{\rho_n}^{\rho_{n+1}}=\arg\min_{T: \Theta \to \Theta} \Big\{ F(T_{\#}\rho_{n})+\dfrac{1}{2\xi}\int_{\Theta} \left\Vert T(\theta)-\theta\right\Vert^2 d\rho_n \Big\}.
\end{equation}
{Notice that the same transport-map parametrization without the $\mathcal{W}_2$ regularization can in principle be used to directly optimize $F$. However, the proximal formulation instead introduces a quadratic Wasserstein regularization that localizes each update around the current iterate. In  particular, since the identity map is feasible, every exact proximal step satisfies $
F(\rho_{n+1})+\frac{1}{2\xi}W_2^2(\rho_{n+1},\rho_n)\leq F(\rho_n)$
so that the displacement of each update is quantitatively controlled by the achieved objective decrease, which allows a quantitative analysis.}
\end{remark}

\subsection{A Wasserstein PL inequality}
\label{subsec:wass_PL}

In this subsection, we introduce the PL inequality in Wasserstein space as in \citep{boufadene2023global}.  
\begin{definition}[Wasserstein Polyak-{\L}ojasiewicz inequality]
\label{Def_PL}
      For any $\rho\in\mathcal{D}(F)$, the objective functional $F$ satisfies the following inequality with $\mu > 0$,
   \begin{equation}
\label{Eqn_wass_PL}
    \int_{\Theta} \left\Vert \nabla \dfrac{\delta F}{\delta \rho}(\rho)\right\Vert^2 d\rho \geq 2\mu (F(\rho)-F^*),
\end{equation}     
where $F^* = F(\rho^*)$ and $\rho^*$ is any global minimizer of $F$. We call~(\ref{Eqn_wass_PL}) a Wasserstein PL inequality.
% Note that the PL inequality implies that any stationary point (not necessarily unique) solving the first-order condition $\nabla \dfrac{\delta F}{\delta \rho}(\rho) = 0, \rho\text{-a.e.},$ is a global minimum $F^*$. 
\end{definition}
The Wasserstein PL inequality generalizes the classical PL inequality in Euclidean space $\left\Vert \nabla f(\theta)\right\Vert^2\geq 2\mu (f(\theta)-f(\theta^*))$ where $f: \Theta\rightarrow \mathbb{R}$~\citep{KarimiNutiniSchmidt_PL}, with a key technical difference that functional on Wasserstein space may not have very good differentiability structure (cf. the discussion after Assumption~\ref{assum_2} below). For KL divergence, the Wasserstein analog of the Euclidean PL inequality is LSI in (\ref{eqn:LSI}), which is almost exclusively used to study the linear convergence for KL divergence-type objective functionals. With a convex function $f$,  the quadratic growth of $f$ implies the PL inequality in Euclidean space~\citep{KarimiNutiniSchmidt_PL} and the quadratic growth of its KL objective implies LSI \citep{yao2024wasserstein}. Previous works show that under certain regularity conditions, the continuous-time dynamics exhibit linear convergence under the Wasserstein PL inequality~\citep{boufadene2023global,kondratyev2016new,chizat2022meanfield}. Our paper considers the problem of minimizing a \emph{general functional} in (\ref{Eqn_objective}), where the convergence analysis of the proximal algorithm is directly based on the Wasserstein PL inequality (\ref{Eqn_wass_PL}). 
\begin{roundedbox}
\begin{assumption}[Minimum selection on proximal trajectory]
\label{assum_2} For every $\rho\in\mathcal{D}(F)$ and every $\rho_\xi \in prox_{F,\xi}(\rho)$, $\nabla\frac{\delta F}{\delta \rho}(\rho_\xi)$ is a subdifferential at $\rho_\xi$ such that
 $$\left\Vert\nabla\frac{\delta F}{\delta \rho}(\rho_\xi)\right\Vert_{L_2(\rho_\xi)}\leq \left\Vert\frac{T_{\rho_\xi}^{\rho}-\mathbf{id}}{\xi}\right\Vert_{L_2(\rho_\xi)}$$
   \end{assumption}
\end{roundedbox}
Some remarks of Assumption~\ref{assum_2} are now in order.

First, it follows from Lemma 10.1.2 in \citep{ambrosio2005gradient} that $\rho_\xi\in D(|\partial F|)$ and $(T_{\rho_\xi}^{\rho}-\mathbf{id}) / \xi$ is a subdifferential of $F$ at $\rho_\xi$ (cf. Definition~\ref{defn:wass_subdifferential} in Supplementary Materials). If $\Theta$ is a compact set in $\mathbb{R}^d$, we have $\nabla\frac{\delta F}{\delta \rho}(\rho_\xi)= ({T_{\rho_\xi}^{\rho}-\mathbf{id}})/{\xi}$ due to the existence of the first variation of $\mathcal{W}_2(\cdot,\rho)$ distance for any fixed $\rho\in\mathcal{P}_2({\Theta})$ (cf. Lemma \ref{lemma_compact}), and thus Assumption \ref{assum_2} automatically holds.  %This property is a key structure in proving a sharper convergence rate for strongly geodesically convex objectives under merely Assumption \ref{assum_1}, see Theorem \ref{Cor_sharp} below.
Moreover, {for non-compact $\Theta$}, both MFLD under the conditions of Corollary \ref{Cor_mfld} and Langevin dynamics under the conditions of Corollary \ref{cor_ld}  satisfy Assumption \ref{assum_2} since $\nabla\frac{\delta F}{\delta \rho}(\rho_\xi)$ is guaranteed to be the subdifferential at $\rho_\xi$ with minimal $L_2(\rho_\xi)$-norm (see proof of Corollary \ref{Cor_mfld} and Corollary \ref{cor_ld}).

Next, we shall highlight throughout this paper that \emph{we do not assume the Wasserstein differentiability, which is an overly restrictive assumption}. For instance, the (geodesically convex) negative entropy functional $\text{Ent}(\rho) = \int \rho \log\rho$ is not Wasserstein differentiable (see Appendix \ref{Appendix_background} for definitions) for {every} $\rho \in \mathcal{P}_2^a(\mathbb{R}^d)$ such that $\text{Ent}(\rho) < \infty$ (cf. Remark 2.26 in~\cite{LBD2025_firstorder}). More precisely, even subdifferential may not exist everywhere and {$\nabla\frac{\delta F}{\delta \rho}(\rho)$ is not guaranteed to be a subdifferential} (e.g., it is not $L_2$ integrable even for $\rho \in \mathcal{P}_2^a(\mathbb{R}^d)$ when $F = \text{Ent}$). Thus, we need to rely on the fact that $F$ is subdifferentiable at $\rho_\xi$ and $(T_{\rho_\xi}^{\rho}-\mathbf{id}) / \xi$ is a subdifferential of $F$ at proximal $\rho_\xi$. {However, we cannot simply relate the subdifferential $(T_{\rho_\xi}^{\rho}-\mathbf{id}) / \xi$ to $\nabla\frac{\delta F}{\delta \rho}(\rho_\xi)$ without Assumption~\ref{assum_2}}, which is a key step to prove our main convergence result in Theorem~\ref{main_thm}. Rather, we connect the proximal iterates $\rho_{\xi}$ (in terms of the subdifferential $(T_{\rho_\xi}^{\rho}-\mathbf{id}) / \xi$) and the PL inequality in (\ref{Eqn_wass_PL}) via the minimal selection principle in Assumption~\ref{assum_2}.

Moreover, unlike the Euclidean space, strong $\mathcal{W}_2$-geodesic convexity generally does not directly imply the Wasserstein PL inequality~(\ref{Eqn_wass_PL}). 
The proof of Theorem \ref{Cor_sharp} doesn't require $\nabla\frac{\delta F}{\delta \rho}(\rho_\xi)$ to be  subdifferential and doesn't rely on the Wasserstein PL inequality. It utilizes the strong geodesic convexity and the subdifferential $(T_{\rho_\xi}^{\rho}-\mathbf{id}) / \xi$ at $\rho_\xi$
 to circumvent relying on Assumption \ref{assum_2}.
\subsection{Convergence rates of exact proximal algorithm}
\label{Subsec_main}

Now, we are ready to state the main theorem of this paper to establish the convergence rate for the Wasserstein proximal algorithm (\ref{Eqn_jko}).

% Theorem~\ref{main_thm} is the one-step contraction rate, which implies the unbiased and linear convergence of the Wasserstein proximal algorithm (\ref{Eqn_jko}).

\begin{theorem}[Convergence rate of the exact proximal algorithm under PL inequality]
\label{main_thm}
Under Assumptions \ref{assum_1} and \ref{assum_2}, if the objective functional $F$ in (\ref{Eqn_objective}) satisfies the PL inequality (\ref{Eqn_wass_PL}), then for any $\xi>0$, the Wasserstein proximal algorithm (\ref{Eqn_jko}) satisfies
\begin{equation}
F(\rho_{n}) - F^*\leq \dfrac{1}{(1+\xi\mu)^{2n}}(F(\rho_0) - F^{*}).
\end{equation}
\end{theorem}

\begin{remark}
It is worth mentioning that our Theorem~\ref{main_thm} extends Theorem 9 in \cite{cyx2022} from the Euclidean space to the Wasserstein space with similar computations. However, we want to point out that the proof of their Theorem 9 is not rigorous. Specifically,~\cite{cyx2022} verified that the Moreau envelope satisfies the Hamilton-Jacobi (HJ) equation through the chain rule on $f_{t,x}(x_t)$ (see definition therein), where the time-derivative of the trajectory $x_t$ is not well-defined because its uniqueness is not guaranteed. In contrast, our Lemma~\ref{Lemma_important} specialized to the Euclidean context can directly lead to satisfaction of the HJ equation, which provides a more rigorous proof than~\cite{cyx2022}.
\end{remark}

Next, we specialize our general-purpose convergence guarantee for the Wasserstein proximal algorithm to MFLD induced from training a two-layer neural network (\ref{eq:twolayer_nn}) in the mean-field regime.
%, and Langevin dynamics, which is a special case of MFLD.
\begin{corollary}[Wasserstein proximal algorithm on MFLD] 
\label{Cor_mfld}
Let $F_\tau$ be the total objective functional in (\ref{eqn:MFLD_obj}). Suppose that the loss function $l(\cdot)$ is either squared loss or logistic loss. If the perceptron $\varphi(\theta; x)$ is bounded by $K$ and $\varphi(\theta;x)$ is smooth in $\theta$, uniformly in $x$,  {then $F_{\tau}$ satisfies a Wasserstein PL inequality with constant $\tau\mu_\tau>0$.} Thus, for any $\xi>0$,
%$F_{\tau}$ satisfies Assumptions \ref{assum_1},  \ref{assum_2} and the PL inequality (\ref{Eqn_wass_PL}). Consequently, 
the Wasserstein proximal algorithm in (\ref{Eqn_jko}) with the objective functional $F_{\tau}$ satisfies
\begin{align}
\label{eqn:mfld_functional}
    F_\tau(\rho_n)-F_\tau^* \leq & \frac{1}{(1+\xi\tau\mu_\tau)^{2n}}(F_\tau(\rho_0)-F_\tau^*), \notag\\
% \label{eqn:mfld_minmizer}
    \mathcal{W}_2(\rho^*,\rho_n)\leq & \sqrt{\frac{2}{\mu_\tau \tau}(F_\tau(\rho_0)-F_\tau^*)} \left( \frac{1}{1+\xi\mu_\tau\tau} \right)^{n}.
\end{align}
\end{corollary}

Our Theorem~\ref{main_thm} can also be applied to derive the convergence rate for backward time-discretized KL divergence (i.e., linear Langevin dynamics).
\begin{corollary}[Wasserstein proximal algorithm on Langevin dynamics]
\label{cor_ld}
Suppose $\nu \propto e^{-f}$ satisfies the $\mu$-LSI condition (\ref{eqn:LSI}) for a potential function $f:\Theta \rightarrow \mathbb{R}$. If $f$ is semiconvex, lower semicontinuous,
% If either of the following holds:
% \begin{enumerate}[noitemsep, topsep=0pt, label=(\alph*)]
%     \item $f$ is semiconvex, lower semicontinuous;
%     % \item $\Theta$ is a compact set and $f$ is lower semicontinuous and bounded from below;
% \end{enumerate}
then for any $\xi>0$,
%Assumptions \ref{assum_1}, \ref{assum_2} and the PL inequality (\ref{Eqn_wass_PL}) are satisfied. Consequently, 
the Wasserstein proximal algorithm in (\ref{Eqn_jko})  with the KL divergence objective functional $\KL$ satisfies 
\begin{align}
\KL(\rho_n \| \nu) \leq & \frac{1}{(1+\mu\xi)^{2n}} \KL(\rho_0 \| \nu), \notag \\
\mathcal{W}_2(\rho_n,\nu) \leq & \sqrt{\frac{2}{\mu} \KL(\rho_0 \| \nu)}\frac{1}{(1+\mu\xi)^{n}}.\notag
\end{align}
\end{corollary}

% \begin{comment}
% \begin{remark}
%     Moreover, if $\Theta$ is a compact set, $f$ is lower semicontinuous and bounded from below, and $\nu=exp(-f)$ satisfies $\mu$-LSI condition (\ref{eqn:LSI}), then Assumptions \ref{assum_1}, \ref{assum_2} and the PL inequality (\ref{Eqn_wass_PL}) are satisfied. Therefore,  for any $\xi>0$, the exact Wasserstein proximal algorithm (\ref{Eqn_jko}) satisfies
% \begin{align}
%     \KL(\rho_n \| \nu) \leq & \frac{1}{(1+\mu\xi)^{2n}} \KL(\rho_0 \| \nu), \notag \\ \mathcal{W}_2(\rho_n,\nu) \leq & \sqrt{\frac{2}{\mu} \KL(\rho_0 \| \nu)}\frac{1}{(1+\mu\xi)^{n}}.\notag
% \end{align}
% \end{remark}
% \end{comment}

    We remark that the objective functions in Corollaries \ref{Cor_mfld} and \ref{cor_ld}
    are linearly convex. In particular, the linear convexity of the MFLD objective
    plays a crucial role in establishing the entropy sandwich argument (cf. proof of Corollary \ref{Cor_mfld}), which in turn
    yields the Wasserstein PL inequality in Corollary \ref{Cor_mfld}.
    This raises a natural question: are there examples satisfying 
    Wasserstein PL inequality that are neither
    linearly convex nor geodesically convex? The following example gives an affirmative
    answer.

\begin{example}[A Wasserstein PL example without linear or geodesic convexity]
\label{Cor_specialexample}
Let \(\Theta=\mathbb R\) and fix \(c\in(0,1)\). Define the functional
\(F:\mathcal P_2(\mathbb R)\to(-\infty,+\infty]\) by
\[
    F(\rho)
    :=
    \begin{cases}
        \displaystyle
        \left(
            m(\rho)+c\sin(m(\rho))
        \right)^2,
        & \rho\in\mathcal P_2^a(\mathbb R), \\[0.8em]
        +\infty,
        & \rho\in\mathcal P_2(\mathbb R)\setminus\mathcal P_2^a(\mathbb R),
    \end{cases}
\]
where
\[
    m(\rho):=\int_{\mathbb R}\theta\,d\rho(\theta).
\]
Then \(F\) satisfies Assumptions \ref{assum_1}, \ref{assum_2}, and  \(F\) satisfies
the Wasserstein PL inequality with constant $\mu=2(1-c)^2$. However, for \(c=1/2\), the functional \(F\) is neither linearly convex nor
geodesically convex.
\end{example}

Next, applying similar proof techniques to the $\mu$-strongly convex objective functional $F$ in (\ref{Eqn_objective}), we obtain a faster convergence rate than those in existing literature, see Remark~\ref{rem:strongconvex_faster} below. In particular, we can avoid evoking Assumption \ref{assum_2} with careful modifications to our proof.

\begin{theorem}[Sharper convergence rates of the exact proximal algorithm: strongly geodesically convex objective]
\label{Cor_sharp}
Under Assumption \ref{assum_1}, if the objective functional $F$ in (\ref{Eqn_objective}) is $\mu$-strongly convex along geodesics, then for any $\xi>0$, the Wasserstein proximal algorithm (\ref{Eqn_jko}) satisfies
\begin{align}
F(\rho_n)-F^* \leq & \frac{1}{(1+\mu\xi)^{2n}}(F(\rho_0)-F^*), \notag \\
\mathcal{W}_2(\rho^*,\rho_n) \leq& \sqrt{\frac{2}{\mu}(F(\rho_0)-F^*)} \frac{1}{(1+\mu\xi)^{n}}. \notag
\end{align}
\end{theorem}
% \begin{remark}
%   \label{rem:strongconvex_proof}  
% \end{remark}
\begin{remark}
\label{rem:strongconvex_faster}
The bound obtained in Theorem~\ref{Cor_sharp} is sharper than those in~\citep{yao2023meanfield} and~\citep{cheng2024} for $\mu$-strongly convex objective, where the convergence rate for functional value $(1+\mu\xi)^{-n} (F(\rho_0)-F^*)$ is proved. Theorem ~\ref{Cor_sharp} applies to functionals corresponding to interacting particle systems $F(\rho)=\int_{\Theta} f(\theta) d\rho(\theta)+\int_{\Theta} w(\theta_1,\theta_2)d\rho(\theta_1)d\rho(\theta_2)+\int_{\Theta} \log\rho(\theta)d\rho(\theta)$, where both the external force $f(\cdot)$ and interaction potential $w(\cdot,\cdot)$ are convex and at least one of them is strongly convex~\citep{yao2024wasserstein}.
\end{remark}

\begin{figure}[t!]
\centering
    \begin{subfigure}[t]{0.4\textwidth}
        \includegraphics[width=\textwidth]{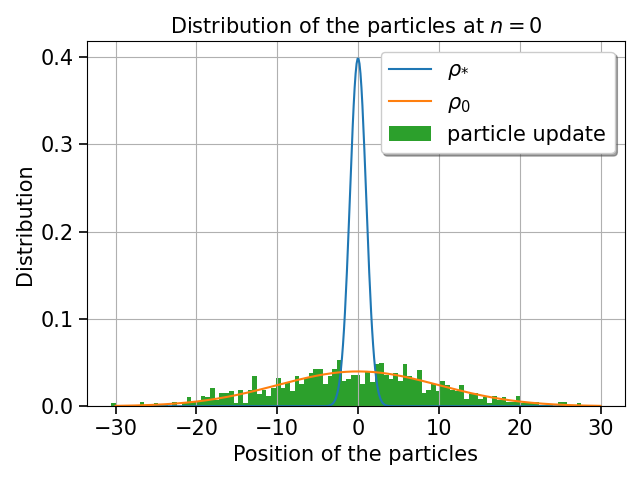}
        \caption{}
    \end{subfigure}
        \begin{subfigure}[t]{0.4\textwidth}
        \includegraphics[width=\textwidth]{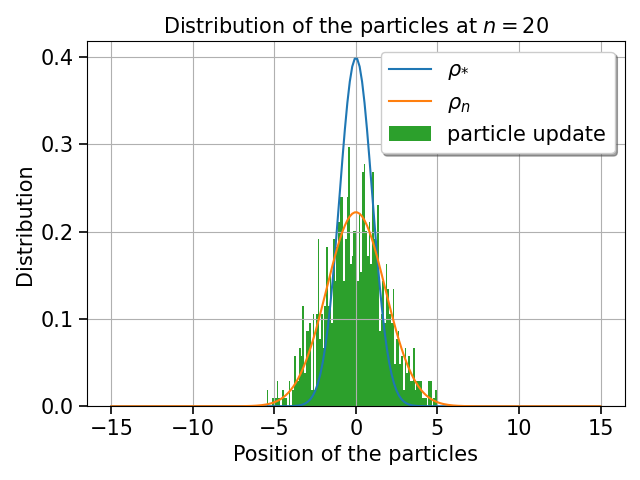}
        \caption{}
        
    \end{subfigure}
            \begin{subfigure}[t]{0.4\textwidth}
        \includegraphics[width=\textwidth]{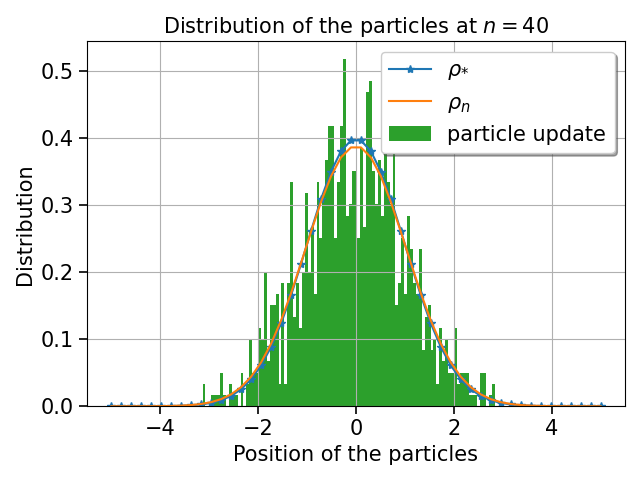}
        \caption{}
    \end{subfigure}
                \begin{subfigure}[t]{0.4\textwidth}
        \includegraphics[width=\textwidth]{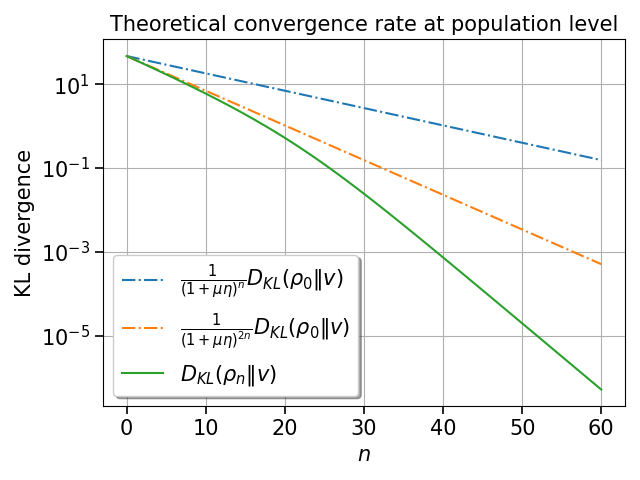}
        \caption{}
        \label{fig_conv_rate}
    \end{subfigure}
    
    \caption{Wasserstein proximal algorithm on the KL divergence with $\nu(\theta) =\exp(-\theta^2/2)$.} 
    % \xc{Is this one-dimensional example?}}
    \vspace{-0.2in}
    \label{Fig_LD}
\end{figure}

\subsection{Convergence rates of inexact proximal algorithm}

Since the trajectory OT map $T_{\rho_n}^{\rho_{n+1}}$ in (\ref{Eqn_prox_op_compute}) is learned through data, this subsection investigates the inexact proximal algorithm where numerical errors are allowed in each iteration. Let ${\overline{\rho}_{n}}$ be the inexact solution of the proximal algorithm at iteration $n$. We need an additional regularity assumption on the proximal flow to provide a quantitative analysis for the proximal algorithm when the {OT map} at each iteration is allowed to be estimated with errors.%on the OT map $T_{\overline{\rho}_n}^{\overline{\rho}_{n+1}}$.
% \textcolor{red}{For instance, if we learn an OT map at each step through a normalizing flow, then Assumption \ref{assum_4} is satisfied, as a normalizing flow is a diffeomorphism. This argument seems not true} \xc{We don't have a proof of this. Probably we can only say empirically this NN or NF can serve as a regularization of the flow.} I provide a "potentially correct" sufficient condition ($T$ is $C^2$ diffeomorphism) in appendix.  
 \begin{roundedbox}
    \begin{assumption}[Regularity of inexact proximal flow]
    \label{assum_4} $\forall n$, $\overline{\rho}_n\in D(|\partial F|)$, $\nabla\dfrac{\delta F}{\delta \rho}(\overline{\rho}_n)\in \partial{F(\overline\rho_n)}$. 

    \end{assumption}
\end{roundedbox}
%%%%%%%%%%%%%%%%%%%%%%%%%%%%%%%%%%%%%%%%%%%%%%%%%%%%%%%
% When the proximal algorithm (\ref{Eqn_jko}) is initialized with ${\rho}_0 \in \mathcal{C}^1(\Theta)$, Assumption~\ref{assum_4} ensures that the inexact proximal flow $\overline{\rho}_n$ remains $\mathcal{C}^1(\Theta)$.  In practice, to optimize over (\ref{Eqn_prox_op_compute}) via $T_{\overline{\rho}_n}^{\overline{\rho}_{n+1}}$, the learned OT map is typically restricted to a specific class, such as a normalizing flow \citep{xu2024normalizing} or a neural network \citep{yao2023meanfield} with a certain structure. 
%%%%%%%%%%%%%%%%%%%%%%%%%%%%%%%%%%%%%%%%%%%%%%%%%%%%%%%%%
{{This assumption is mainly for technical purposes, which guarantees that $\nabla\frac{\delta F}{\delta \rho}(\overline{\rho}_n)$ is a subdifferential}. Next we define the error $\beta_{n+1}:\Theta\rightarrow \Theta$ in estimating the Wasserstein subdifferential in the $(n+1)$-th iterate as,
\begin{equation}
\label{Eqn_inexact}
    \beta_{n+1}=\xi^{-1}(T^{\overline\rho_n}_{\overline\rho_{n+1}}-\mathbf{id})- \nabla \dfrac{\delta F}{\delta \rho}(\overline\rho_{n+1}).
\end{equation}
% under reasonable regularity such as $\rho\in \mathcal{C}^1(\Theta)$ (cf. Lemma 10.4.1 in \cite{ambrosio2005gradient}).
The definition~(\ref{Eqn_inexact}) of inexact subdifferential follows from similar ideas as in~\citep{cheng2024,yao2023meanfield}.{ When $\Theta$ is compact, if ${\overline\rho_{n+1}}$ is the exact solution of the Wasserstein proximal algorithm (\ref{Eqn_jko}), then $\beta_{n+1}=0$ by previous discussion. For a general non-compact domain, such an identification is not automatic,
since the subdifferential need not be unique. However, it is still reasonable to use $\beta_{n+1}$ to characterize the inexactness of proximal iteration \citep{cheng2024}}, and we use
\begin{equation}
\int_{\Theta} \Vert\beta_{n+1}\Vert^2d\overline\rho_{n+1}\leq \epsilon_{n+1},
\end{equation} as an error criterion for the inexact proximal update. Furthermore, we need the additional geodesic semiconvexity assumption apart from PL inequality to control the inexact error, which differs
from Section \ref{Subsec_main} that only relies on PL inequality. {Examples that satisfy both semiconvexity
and PL inequality include the objective of MFLD under the assumptions of Corollary \ref{Cor_mfld} and KL divergence objective under the assumptions of Corollary \ref{cor_ld}}. Now we are in the position to quantify the impact of numerical errors on the proximal algorithm, which demonstrates how the inexact error impacts the convergence behavior under PL inequality. 
\begin{theorem}[Inexact proximal algorithm under PL inequality]
\label{Thm_inexact}
Under Assumptions \ref{assum_1} and \ref{assum_4}, if $F$ is $(-L)$- geodesically semiconvex, $F$ satisfies the PL inequality (\ref{Eqn_wass_PL}), $\epsilon_{n}\leq \epsilon$ for $n\in\mathbb{N}$ and $0 < \xi \leq \frac{1}{L}$, then we have,\\
\begin{equation}
     F(\overline\rho_n)-F^*\leq\dfrac{1}{(1+\mu\xi)^n}(F(\rho_0)-F^*)+\dfrac{\epsilon}{2\mu}.
\end{equation}
\end{theorem}

We remark that the numerical error bound $\epsilon$ does not have to be small and does not need to converge to zero. We shall explore the behavior of the numeric error through simulations in Section~\ref{subsec_mfld_exp}. 

To make our complexity analysis complete, we give algorithmic convergence rates when we have increasing accuracy to solve the proximal subproblems as the algorithm progresses (i.e., $\epsilon_n \to 0$ as $n \to \infty$). 

% Such a regime has been empirically verified with decreasing $\epsilon_n$ in the JKO-training of neural networks~\cite{cheng2024}.

\begin{corollary}[Inexact proximal algorithm under PL inequality with increasing accuracy]
\label{Cor_inexact}
Under the assumptions of Theorem~\ref{Thm_inexact}, we have the following.
     (a) If $\epsilon_n\leq C_{exp}\gamma^n$ with $\gamma\in(0,1)$, then there exists $C_1=C_1(\mu,\xi,\gamma,C_{exp})$
  such that % \dfrac{BC_{exp}}{\vert A-\xi\vert}
 \begin{equation}
 \begin{aligned}
     F(\overline\rho_n)-F^*\leq & \dfrac{1}{(1+\mu\xi)^n}(F(\rho_0)-F^*)
     +C_1 \max\left\{\frac{1}{1+\mu\xi},\dfrac{1+\gamma}{2} \right\}^{n+1}.
     \end{aligned}
 \end{equation}
 (b) If $\epsilon_n\leq C_{poly}n^{-\zeta}$ with $\zeta>0$, then there exists $C_2=C_2(\mu,\xi,\zeta,C_{poly})$ such that
 \begin{equation}
     F(\overline\rho_n)-F^*\leq \dfrac{1}{(1+\mu\xi)^n}(F(\rho_0)-F^*)+\dfrac{C_2}{n^{\zeta}}.
 \end{equation}
\end{corollary}
Corollary \ref{Cor_inexact} gives a quantitative decay of how the inexact error impacts the convergence behavior under the PL inequality. 

\begin{figure*}[t!]
\centering
    \begin{subfigure}[t]{0.41\textwidth}
        \includegraphics[width=\textwidth]{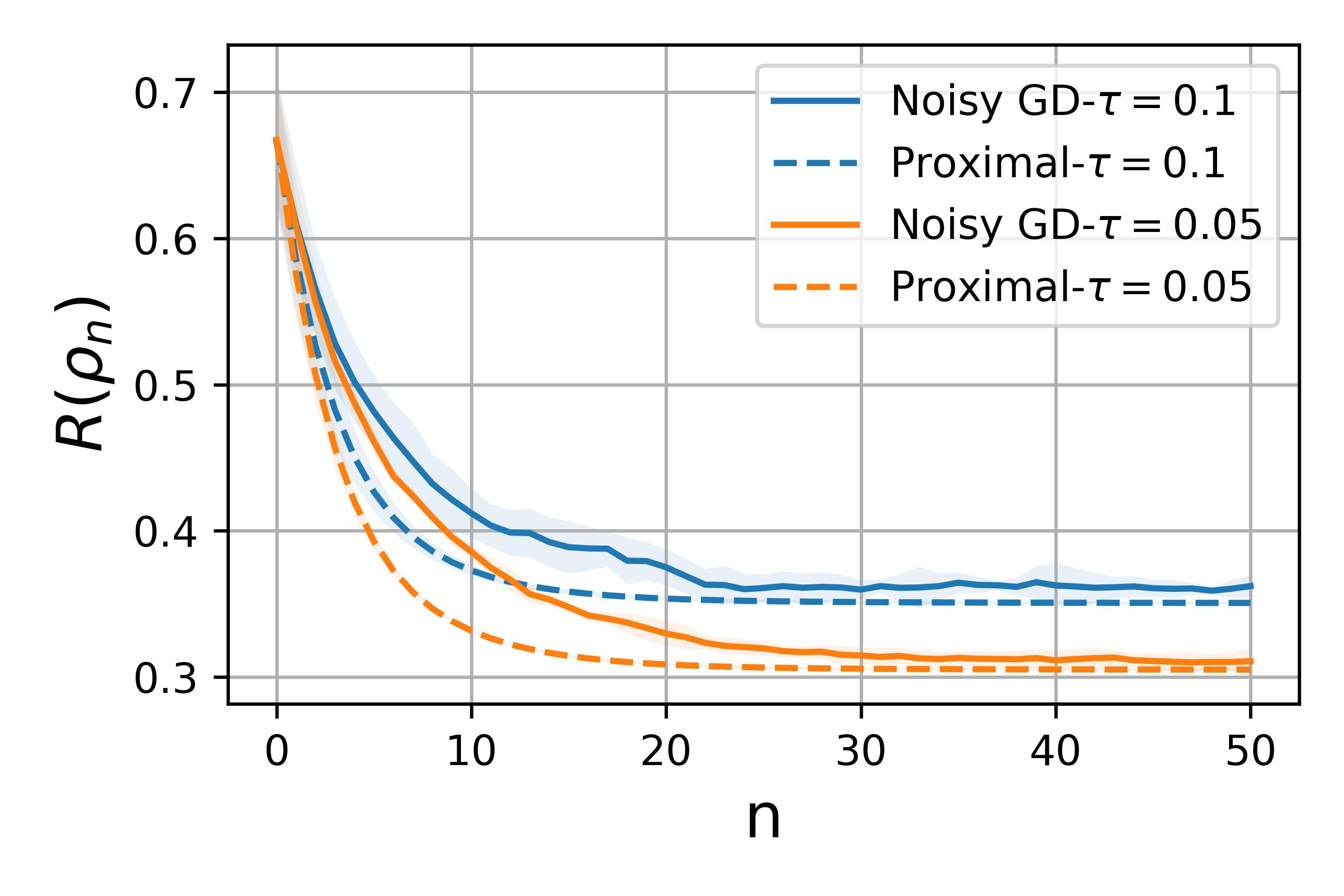}
        \caption{}
        \label{Fig_loss}
    \end{subfigure}
        \begin{subfigure}[t]{0.41\textwidth}
        \includegraphics[width=\textwidth]{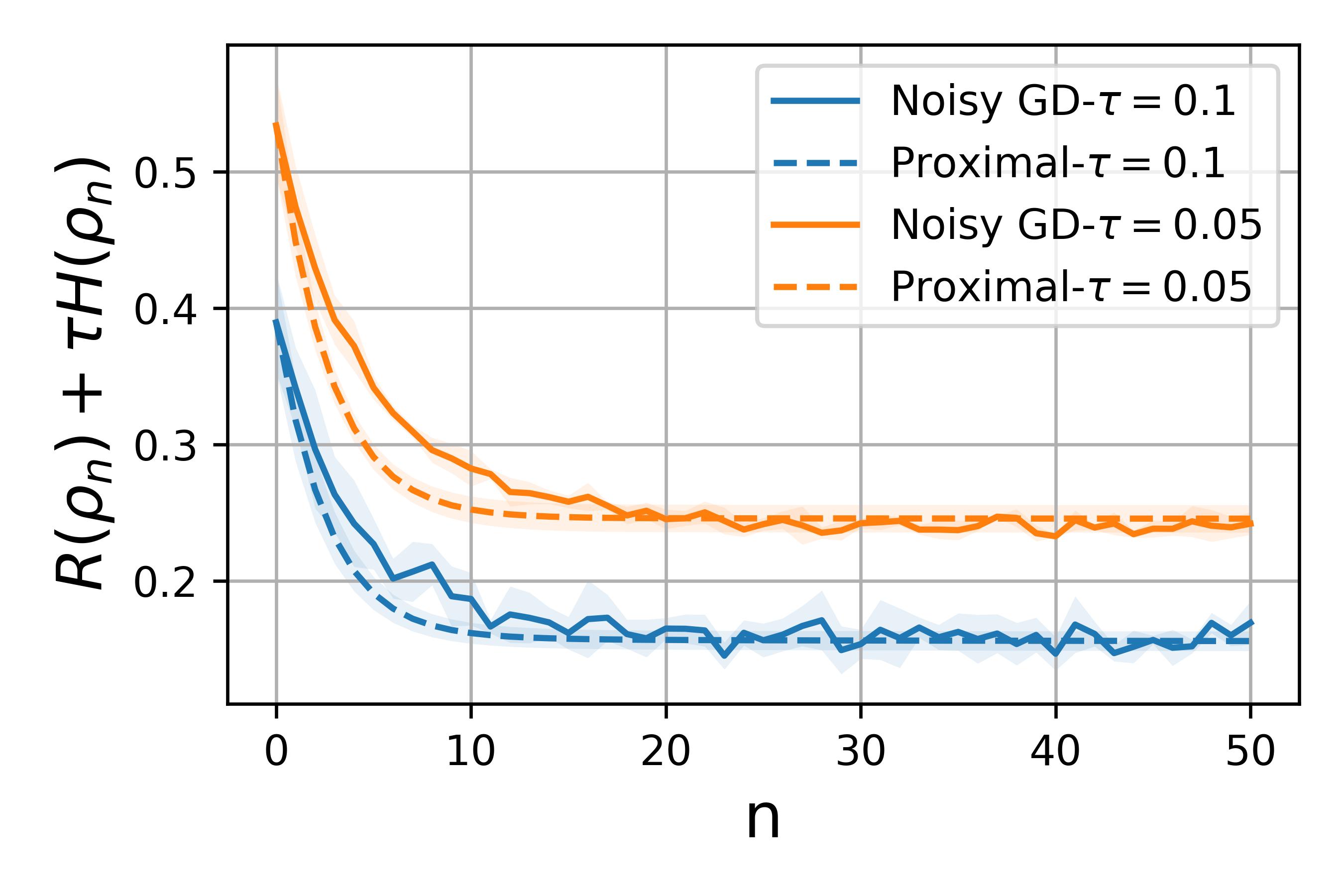}
        \caption{}
        \label{Fig_total_loss}
    \end{subfigure}
            \begin{subfigure}[t]{0.41\textwidth}
        \includegraphics[width=\textwidth]{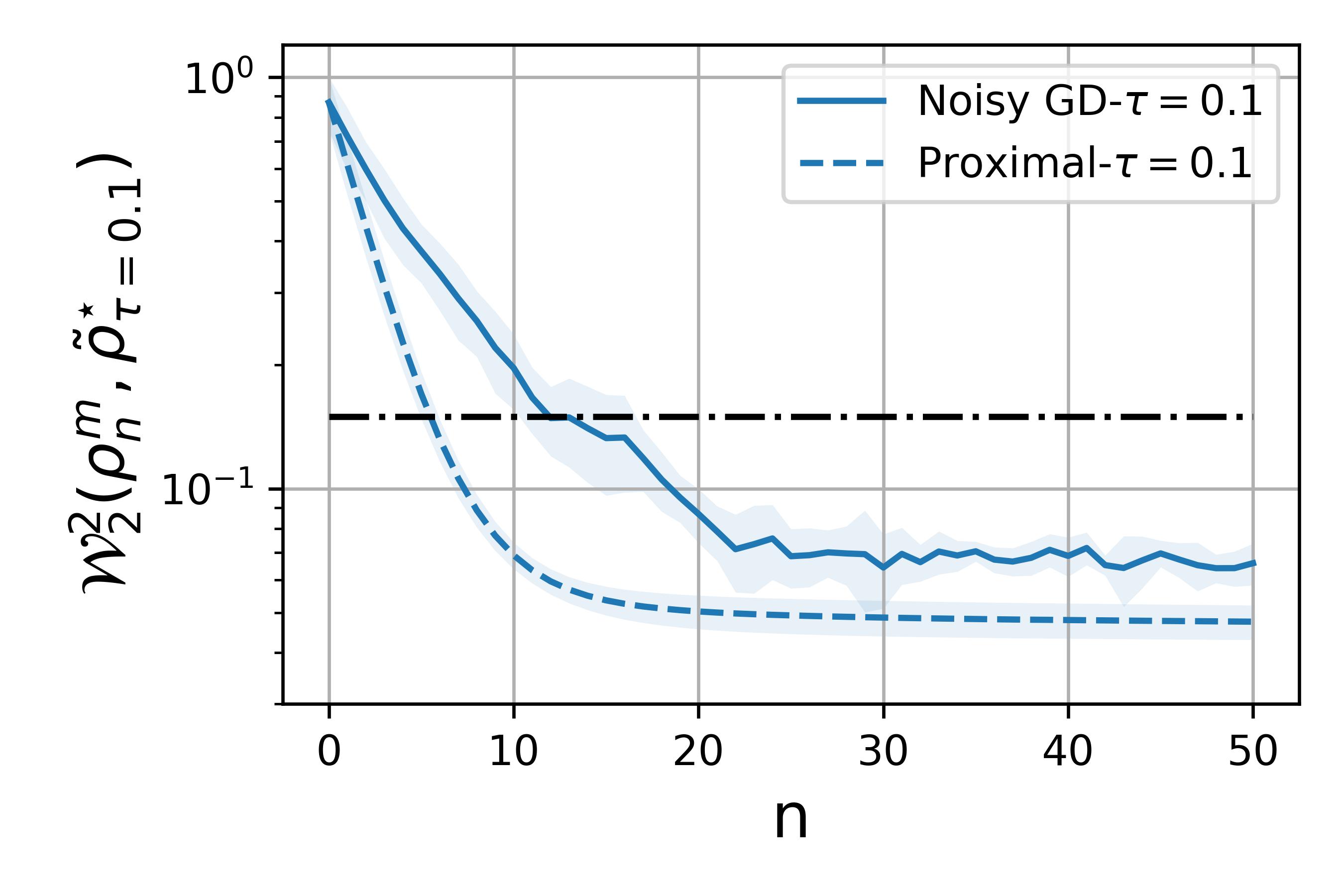}
        \caption{}
        \label{Fig_W2_0.1}
    \end{subfigure}
                \begin{subfigure}[t]{0.41\textwidth}
        \includegraphics[width=\textwidth]{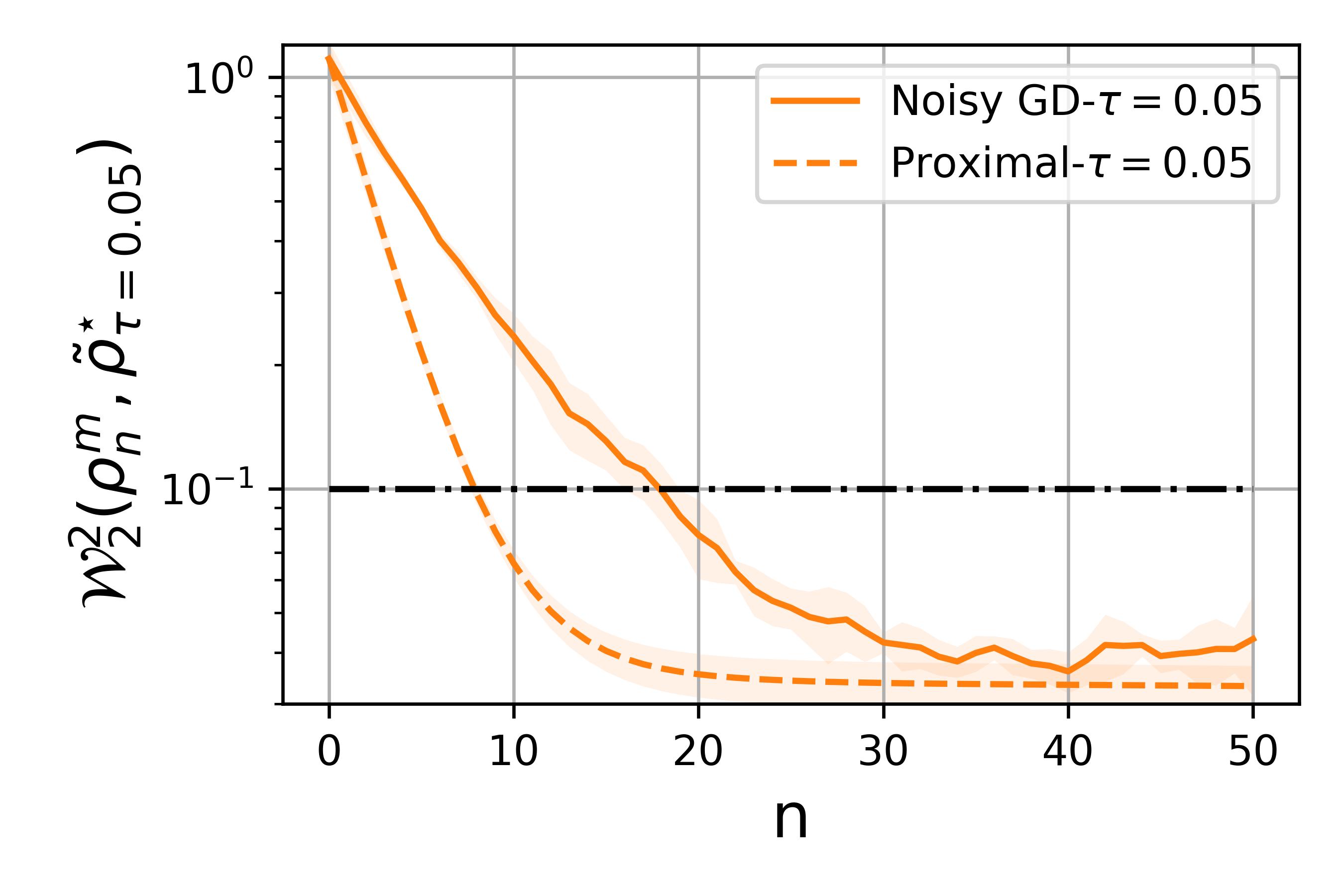}
        \caption{}
        \label{Fig_W2_0.04}
    \end{subfigure}
    
    \caption{Wasserstein proximal algorithm on MFLD objective. Note that for Figure \ref{Fig_W2_0.1} and Figure \ref{Fig_W2_0.04}, the y-axis is on log scale.} \vspace{-0.2in}
    \label{Fig_MFLD}
\end{figure*}
\section{Numerical experiments}
\label{Sec_apply}

In this section, we present two applications of the proximal algorithm. The first is on the KL divergence functional, where the particle and the distribution updates can be computed explicitly. The second is on the proximal algorithm training for the regularized objective of two-layer mean-field neural networks.

\subsection{Linear Langevin dynamics}
We apply the proximal algorithm on KL divergence with the target distribution $\nu(\theta) \propto\exp(-\theta^2 / 2)$ where $\theta\in \mathbb{R}$. Note that $D_{KL}(\cdot\Vert \nu)$ is $1$-strongly convex along geodesics. %We provide numerical experiments to explore the dynamical behavior of the proximal algorithm.
When both initialization and the target distributions are Gaussian, problem (\ref{Eqn_prox_op_compute}) can be explicitly solved and $\rho_n$ remains Gaussian for every $n$. Closed forms of the particle and distribution updates are available in \citep{wibisono18a}.  Additionally, the $\mathcal{W}_2$ distance between two Gaussian distributions can be computed explicitly.  In the experiment, we set the initialization Gaussian distribution to be $\mathcal{N}(0,100)$, step size $\xi=0.1$, and iterations equal to 60. Our implementation is based on the implementation of \cite{salim2020wasserstein}. 

In Figure \ref{Fig_LD}, the distribution of particles, represented by the histogram, approximates $\rho_n$ and converges to $\rho^*$ after several iterations ($\approx$40 iterations in this experiment). The linear convergence result of $D_{KL}(\rho_n\Vert \nu)$ in Figure \ref{fig_conv_rate} demonstrates a sharper bound holds for $\mu$-convex objective with respect to \citep{yao2023meanfield,cheng2024}, as Theorem \ref{Cor_sharp} suggests.

\subsection{Mean-field neural network training with entropy regularization} 
\label{subsec_mfld_exp}
In practice, the optimization problem (\ref{Eqn_prox_op_compute}) typically lacks an explicit solution. {Therefore, in contrast to our population-level theory, which is formulated over
absolutely continuous probability measures, the numerical implementation
uses a finite-particle approximation }of the time-evolving probability distributions, and we can solve an approximate $T^{n+1}$ using functional approximation methods (cf. more details in Appendix Section~\ref{subsec:particle+functional_app}). 
In this experiment, we minimize the MFLD entropy-regularized total objective (\ref{eqn:MFLD_obj}) with $\varphi(\theta,x)=\tanh(\theta^T x)$, which satisfies the Wasserstein PL inequality (see Corollary \ref{Cor_mfld}), but is not geodesically convex. The parameters are set as $d=2$, $\lambda = 0.1$, $\tau = \{0.05, 0.1\}$, with the number of particles $m = 100$ and a discretized step size of $\xi = 0.2$.  We generate $N=1000$ training data samples using a teacher model $y=\sin(\alpha^Tx)$, where $x\sim \mathcal{N}(0,I)$. 

Our goal is to compare the proximal algorithm with the neural network-based functional approximation (\ref{Eqn_learn_otmap}) with the noisy gradient descent algorithm (\ref{noise_gd_eqn}). We first randomly generated a dataset and then conducted 5 repeated experiments for both $\tau = 0.05$ and $\tau = 0.1$ on the same generated data. For each 
experiment, a new weight empirical distribution is generated from the standard Gaussian distribution, and both algorithms will use the same weight as the initial value. To learn the optimal transport map, we train a RealNVP normalizing flow  for each step using SGD optimizer with learning rate 0.005 and 150 iterations by default unless explicitly mentioned (see Appendix \ref{app:further_exp} for more details).

\begin{figure}
    \centering
    \includegraphics[width=0.42\textwidth]{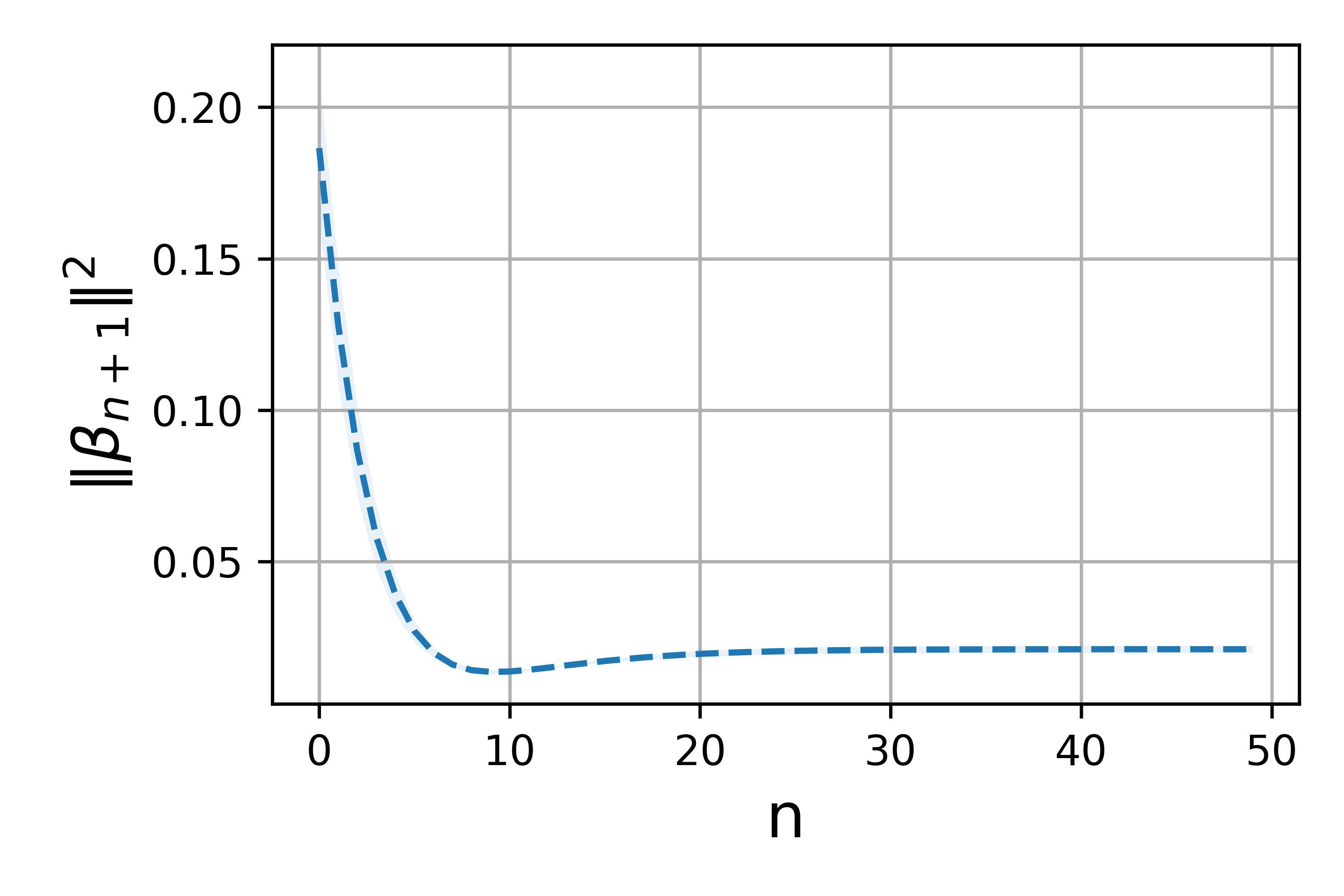}
    \vspace{-6mm}
    \caption{Behavior of inexact error $\Vert\beta_{n+1}\Vert^2_{{L_2(\rho_{n+1}^m)}}$ ($\tau=0.1, \xi=0.2$).}
    \label{Fig_inexact}
    \vspace{-0.2in}
\end{figure}

In Figure \ref{Fig_loss} and Figure \ref{Fig_total_loss}, we observe that both the $L^2$-regularized loss $R$ and the total objective $F_{\tau}$ converge under two algorithms, where the nearest neighbor estimator \citep{kozachenko1987sample} is used to estimate $\int_{\Theta} \rho \log\rho$. To better depict the convergence rate of the Wasserstein proximal algorithm and the Langevin algorithm, we obtain a reference $\tilde{\rho}_\tau^*$ of $\rho^*$,  by running the noisy gradient descent algorithm with very small step size $0.001$ and $m=1000$ particles.  
In the early training phase of both algorithms, $\mathcal{W}_2^2(\rho_n^m,\tilde{\rho}_\tau^*)$ is dominated by $\mathcal{W}_2^2(\rho_n^m,{\rho}^*)$ and exhibits a linear convergence above the black dash-dot line as shown in Figure \ref{Fig_W2_0.1} and Figure \ref{Fig_W2_0.04}. Within this phase, the Wasserstein proximal algorithm demonstrates a faster linear rate thanks to the unbiased linear convergence nature of $\mathcal{W}_2^2(\rho_n^m,\rho^*)$ (cf. Corollary \ref{Cor_mfld}). However, both algorithms have similar biases at convergence, for which we conjecture that the particle discretization error of $\rho_n^m$ dominates $\mathcal{W}_2^2(\rho_n^m,\tilde{\rho}_\tau^*)$ while close to convergence. We validate our
conjecture through further experiments in Appendix \ref{app:further_exp}. {We emphasize that our numerical comparison primarily concerns convergence with respect to the number of iterations rather than wall-clock time. In our current implementation, each proximal iteration takes around $200\times$ the runtime of one noisy-GD step because it requires training a transport map. However, the proximal algorithm has a practical advantage: the learned optimal transport maps can be saved and reused without retraining \cite{xu2024normalizing}, a feature that is particularly useful in generative modeling and sampling tasks.} 

We also explore the empirical behavior of $\Vert\beta_{n+1}\Vert^2_{L_2(\rho_{n+1}^m)}$ as $n$ increases in Figure \ref{Fig_inexact}. $\nabla \log \rho$ is approximately computed by a kernel approach, where we use a Gaussian kernel
with a bandwidth 0.5. We observe that the approximate $\Vert\beta_{n+1}\Vert^2_{L_2(\rho_{n+1}^m)}$ is bounded, suggesting the boundedness assumption $\epsilon_n \leq \epsilon$ in Theorem \ref{Thm_inexact}.

\section{Conclusion}
In this work, we provide a convergence analysis of the Wasserstein proximal algorithm under a Wasserstein PL inequality, without assuming geodesic convexity. Our analysis also improves upon existing convergence rates when strong geodesic convexity holds. We further analyze inexact gradient variants under an additional geodesic semiconvexity condition. When applied to proximal training of mean-field neural networks, we prove linear convergence of the entropy-regularized total objective, and show empirically that proximal training converges faster than the noisy gradient descent algorithm.

\subsubsection*{Acknowledgments}
We would like to thank the anonymous reviewers for their constructive and insightful comments. We also acknowledge the use of ChatGPT in constructing Example 3.9 and in proofreading the manuscript. This work was partially supported by NSF DMS-2413404 and a gift from the Simons Foundation.
\bibliography{main}
\bibliographystyle{tmlr}

%%%%%%%%%%%%%%%%%%%%%%%%%%%%%%%%%%%%%%%%%%%%%%%%%%%%%%%%%%%%%%%%%%%%%%%%%%%%%%%
%%%%%%%%%%%%%%%%%%%%%%%%%%%%%%%%%%%%%%%%%%%%%%%%%%%%%%%%%%%%%%%%%%%%%%%%%%%%%%%
% APPENDIX
%%%%%%%%%%%%%%%%%%%%%%%%%%%%%%%%%%%%%%%%%%%%%%%%%%%%%%%%%%%%%%%%%%%%%%%%%%%%%%%
%%%%%%%%%%%%%%%%%%%%%%%%%%%%%%%%%%%%%%%%%%%%%%%%%%%%%%%%%%%%%%%%%%%%%%%%%%%%%%%
\newpage
\appendix
\onecolumn
\section{Background on optimal transport and Wasserstein space}
\label{Appendix_background}
\subsection{Wasserstein distance and optimal transport}
The squared 2-Wasserstein distance is defined as the solution to the Kantorovich problem
\begin{equation}
\label{Eqn_ot_kt}
    \mathcal{W}_2^2(\rho,\tilde\rho) := \min_{\pi\in \Pi(\rho,\tilde\rho)}\int_{\Theta \times \Theta} \left\Vert \theta-\tilde\theta\right\Vert^2 d\pi(\theta,\tilde \theta),
\end{equation} where  $\Pi(\rho,\tilde\rho)\subset \mathcal{P}_2(\Theta\times \Theta)$ is the set of all coupling distributions with marginals $\rho$ and $\tilde\rho$. The optimal solution $\pi^*$ is called the optimal transport plan. When $\rho\in\mathcal{P}_2^a(\Theta)$, it is known from celebrated Brenier's theorem that the solution $T_\rho^{\tilde\rho}$ of Monge's problem exists, and the optimal transport plan is $\pi^* = (\mathbf{id}, T_\rho^{\tilde{\rho}})_\#\rho$.

\subsection{Wasserstein subdifferential and $\mu$-convex functionals}
% Throughout this paper, we assume that 
% $\mathcal D(F)\subset\mathcal{P}_2^a(\Theta)$ for simplicity in this paper, which simplifies the following definitions. We will revisit this assumption in Section \ref{Sec_main}.
\begin{definition}[Fr\'echet subdifferential, Definition 10.1.1 in \cite{ambrosio2005gradient}]\label{defn:wass_subdifferential}
Let $F:\mathcal{P}_2(\Theta)\rightarrow(-\infty,+\infty]$ be proper, lower semicontinuous, and let $\rho\in D(|\partial F|)$ where $|\partial F|(\rho)$ denotes the metric slope \citep{ambrosio2005gradient}. We say that $\mathbf{v}\in L^2(\rho;\Theta)$ belongs to the (Fr\'echet) subdifferential $\partial F(\rho)$ of $F$, written as $\mathbf{v} \in \partial F(\rho)$, if
$$F(\tilde{\rho})\geq F(\rho)+\int_{\Theta}\langle \mathbf{v}(\theta),T_{\rho}^{\tilde{\rho}}(\theta)-\theta \rangle d\rho(\theta)+o(\mathcal{W}_2(\rho,\tilde{\rho})).$$  
 Also, $\mathbf{v}\in \partial F(\rho)$ with minimal $\left\Vert\cdot \right\Vert_{L_2(\rho)}$ norm is denoted as $\partial^{\circ}F(\rho)$, and it is unique (cf. Lemma 10.1.5 in \cite{ambrosio2005gradient}). If $\partial F(\rho)$ is not empty, then we say $F$ is subdifferentiable at $\rho$. $F$ is differentiable at $\rho$ if both $F$ and $-F$ are subdifferentiable.

\end{definition}
\begin{definition}[First variation] 
\label{def:first_variation}
Let $F:\mathcal{P}_2(\Theta)\rightarrow(-\infty,+\infty]$ be proper, lower semicontinuous, and $\rho\in D(F)$. The first variation $\frac{\delta F}{\delta \rho}(\rho):\Theta\rightarrow \mathbb{R}$ is defined as
\begin{equation}
    {\lim_{\epsilon\rightarrow 0^{+}}}\dfrac{F({\rho}+\varepsilon\chi)-F(\rho)}{\varepsilon}=\lim_{\varepsilon\rightarrow 0^+}
\frac{
F((1-\varepsilon)\rho+\varepsilon\widetilde\rho)-F(\rho)
}{\varepsilon}=\int\dfrac{\delta F}{\delta \rho}(\rho)d\chi
\end{equation}
for any perturbation $\chi=\tilde{\rho}-\rho$ such that  $\tilde{\rho}\in\mathcal{P}_2(\Theta)$ and $(1-\varepsilon)\rho+\varepsilon \tilde{\rho}\in D(F)$.
\end{definition}

% Note that if $\xi\in\partial F$ for some smooth functional $F$, then under mild conditions, e.g., Lemma 10.4.1 in~\cite{ambrosio2005gradient}, we have with $\rho$-probability one,
% \begin{equation}
% \xi(\theta) = \nabla\dfrac{\delta{F}}{\delta \rho}(\rho)(\theta),
% \end{equation}
% which implies that the Wasserstein gradient of $F$ at $\rho$ is $\nabla\frac{\delta F}{\delta\rho}(\rho)$ . 

\begin{definition}[$\mu$-convexity along geodesic, {Section 10.1.1 of \cite{ambrosio2005gradient}}] 
\label{def:mu_conxity}
A proper, lower semicontinuous functional $F$ is called $\mu$-convex along geodesic ($\mu\in\mathbb{R}$) if there exists a constant-speed Wasserstein geodesic
$(\rho_t){t\in[0,1]}$ connecting $\rho_0,\rho_1\in D(F)$  such that
\begin{equation*}
\label{eqn}
F(\rho_t)
\leq
(1-t)F(\rho_0)+tF(\rho_1)
-\frac{\mu}{2}t(1-t)\mathcal{W}_2^2(\rho_0,\rho_1),
\qquad t\in[0,1].
\end{equation*} If $F$ is $\mu$-convex along geodesics, then at $ \rho\in\mathcal{D}(F)$,   
\begin{equation}
\label{eqn:mu_conxity}
 F(\tilde{\rho})\geq F(\rho)+\int_{\Theta}\langle \mathbf{v}(\theta),T_{\rho}^{\tilde{\rho}}(\theta)-\theta \rangle d\rho(\theta)+\dfrac{\mu}{2}{\mathcal{W}_2^2(\rho,\tilde{\rho})},\forall \tilde{\rho}\in \mathcal{D}(F)
\end{equation}
if and only if $\mathbf{v}\in\partial F(\rho)$. In particular, if $\mu\geq0$, we call F \textbf{geodesically convex}; If $\mu<0$, we call $F$ \textbf{geodesically semiconvex}. We refer to Section 9 of \cite{ambrosio2005gradient} for definitions of $\mu$-convexity (and semiconvexity) in Euclidean space, which are similar to the definition above.
\end{definition}

\begin{definition}[Weak convergence]Let $\mathcal{C}_b(\Theta)$ be the set of all continuous bounded functions on $\Theta$ and  $\mathcal M(\Theta)$ be the set of all finite signed measures on $\Theta$. We say that $\rho_k\in \mathcal{M}(\Theta)$ converges to $\rho\in\mathcal{M}(\Theta)$ weakly if for every $\phi\in \mathcal{C}_b(\Theta)$,
\begin{equation}
    \lim_{k\rightarrow \infty} \int_\Theta \phi d\rho_k = \int_\Theta \phi d\rho.
\end{equation}
The weak convergence is also called narrow convergence in the literature \cite{santambrogio2015optimal}.
\end{definition}

\section{Technical lemmas}
\label{app:techLemmas}

% \begin{lemma}[An equation connecting $\partial_\xi u(\theta,\xi)$ with $\mathcal{W}_2^2( \rho_\xi,\rho)$]
% If $F(\rho) < \infty$, the following equation holds for $\xi\in(0,+\infty)$ with at most countable exceptions,
% \label{Lemma_important}
% \begin{align}
%     \partial_{\xi}u(\rho,\xi)&=-\frac{1}{2\xi^2}\mathcal{W}_2^2(\rho_\xi,\rho). \label{Eqn_important_lemma}
% \end{align}
% \end{lemma}
% Lemma \ref{Lemma_important} relates the time-derivative of the Hopf-Lax semigroup to the squared Wasserstein distance between proximal and initial points. Proof of Lemma \ref{Lemma_important} can be found in Appendix~\ref{app:proofs}, which essentially follows from \cite{Ambrosio_2013}.
% \begin{lemma}
% \label{Lemma_opt}
% For the Wasserstein proximal operator with $\xi > 0$,
% \begin{equation}
%     \rho_{\xi}=prox_{F,\xi}(\rho) = \argmin_{\tilde{\rho}\in\mathcal{P}_2({\Theta})}F(\tilde{\rho})+\dfrac{1}{2\xi}\mathcal{W}_2^2(\tilde{\rho},\rho),
% \end{equation}
% the first-order optimality condition is given by
% \begin{equation}
% \label{Eqn_optoptfirst}
% \nabla \dfrac{\delta F}{\delta \rho}(\rho_{\xi})=\dfrac{T_{\rho_{\xi}}^{\rho}-\mathbf{id}}{\xi}.
% \end{equation}
% \end{lemma}

% \begin{proof}[Proof of Lemma~\ref{Lemma_opt}]As $\mathcal{D}(F)\subset \mathcal{P}_2^a(\Theta)$, then $\rho_\xi\in \mathcal{P}_2^a(\Theta)$. Lemma 10.1.2 in~\cite{ambrosio2005gradient} implies that $\dfrac{T_{\rho_{\xi}}^{\rho}-\mathbf{id}}{\xi}\in \partial F(\rho_\xi)$ and thus (\ref{Eqn_optoptfirst}) holds $\rho_\xi$-almost everywhere.
% \end{proof}

Recall that for a general metric space $(Z,d)$, the Hopf-Lax semigroup is defined as
  \begin{equation}
      u(z,\xi)=\inf_{z'\in Z} \left\{ f(z')+\dfrac{1}{2\xi}d(z',z)^2 \right\}
  \end{equation} 
for  $z\in Z$ such that $f(z)<+\infty$.
  
\begin{lemma}[Proposition 3.1 and 3.3 in~\cite{Ambrosio_2013}]
\label{lem:ambrosio_tech}
  Let
  \begin{equation}
      D_{+}(z,\xi) := \sup\mathop{\lim\sup}_{n\rightarrow \infty}d(z,z_n' ), \qquad D_{-}(z,\xi) := \inf\mathop{\lim\inf}_{n\rightarrow \infty}d(z,z_n' ),
  \end{equation}
  where the supremum and the infimum run among all minimizing sequences $(z_n')$.  For $z\in Z$ such that $f(z)<+\infty$, we define $t^*(z)=\sup\{t>0: u(z,t)>-\infty\}$. Then if $f(z)<+\infty$, we have
  
  (a) $D_{+}(z,\xi)=D_{-}(z,\xi)$ holds for $\xi\in (0,t^*(z))$ except for at most countable exceptions;
  
(b)  If and only if $D_{+}(z,\xi)=D_{-}(z,\xi)$, the map $\xi\rightarrow u(z,\xi)$ is differentiable in $(0,t^*(z))$  and 
$$\dfrac{d}{d\xi}u(z,\xi)=-\dfrac{D^2_{+}(z,\xi)}{2\xi^2}=-\dfrac{D^2_{-}(z,\xi)}{2\xi^2}.$$
\end{lemma}

\begin{lemma}[Existence of minimizer]
\label{Lemma_exist} 
If $F:\mathcal{P}_2(\Theta)\rightarrow (-\infty,\infty]$ is proper, weakly lower semicontinuous {and lower bounded}, then proximal algorithm (\ref{Eqn_jko}) admits a minimizer.
\begin{proof}
{Without loss of generality, let $(\widetilde{\rho}_k)_{k\in\mathbb N}$ be a minimizing sequence such that $\widetilde{\rho}_0=\rho$ and
\[
F(\widetilde{\rho}_{k+1})
+\frac{1}{2\xi}W_2^2(\widetilde{\rho}_{k+1},\rho)
\leq
F(\widetilde{\rho}_k)
+\frac{1}{2\xi}W_2^2(\widetilde{\rho}_k,\rho)
\]
for every $k\in\mathbb N$. Since \(F\) is bounded from below, say \(F\geq M\), for all sufficiently large \(k\),
\[
    M+\frac{1}{2\xi}\mathcal W_2^2(\widetilde\rho_k,\rho)
    \leq
    F(\widetilde{\rho}_k)+\frac{1}{2\xi}\mathcal W_2^2(\widetilde\rho_k,\rho)\leq 
    F(\rho)
\]
Therefore, \((\widetilde\rho_k)_{k\in\mathbb N}\) is contained in
\[
    \mathcal{B}_r(\rho)
    :=
    \{\nu\in\mathcal P_2(\Theta)\mid \mathcal W_2^2(\rho,\nu)<r\}
\]
for $r=2\xi(F(\rho)-M)+1$.} By the weak lower semicontinuity of both $F$ and
$\mathcal{W}_2^2(\cdot,\rho)$, every weak limit point of the
minimizing sequence is a minimizer of the proximal problem. Thus, to obtain a weakly convergent subsequence of the minimizing
sequence, it remains to show that \(\mathcal B_r(\rho)\) is weakly precompact. 
  For the sake of completeness, we provide a proof here. The proof is essentially contained in Section 10.1 of \cite{ambrosio2005gradient}. By Prokhorov's theorem, it suffices to prove the tightness of $\mathcal{B}_r(\rho)$, i.e.,  there is a sequence of compact sets $(K_i)_{i\in\mathbb{N}}$ such that
\begin{equation}
    \nu(\Theta \backslash K_i )\leq 1/i, \quad \forall \nu\in \mathcal{B}_r(\rho).
\end{equation}
 We prove for any $\varepsilon>0$, we can find compact set $K$ such that 
\begin{equation}
    \nu(\Theta \backslash K )\leq \varepsilon, \quad \forall \nu\in \mathcal{B}_r(\rho).
\end{equation}   
by contradiction. Assume there exists $\varepsilon$, for any compact $K$, there exists $\nu\in \mathcal{B}_r(\rho)$ such that $\nu(\Theta\backslash K)>\varepsilon$.  As a singleton $\{\rho\}$ constitutes a tight family, we can find a compact set $K_{\rho,\varepsilon}$ such that 
$$\rho(\Theta \backslash K_{\rho,\varepsilon})<\varepsilon/2.$$
We define a compact set $U_R=\{\tilde\theta| \min_{\theta\in K_{\rho,\varepsilon}}\left\Vert\theta-\tilde\theta\right\Vert^2\leq R, \}$. Then there exists $\nu$ such that $\nu(\Theta\backslash U_{3r/\varepsilon})> \varepsilon$. However, $\mathcal{W}_2^2(\rho,\nu)\geq \frac{\varepsilon}{2}\frac{3r}{\varepsilon}=\frac{3}{2}r$, contradiction. 
% Furthermore, note that since $\Theta$ is not compact, $\mathcal{B}_r(\rho)$ under the Wasserstein metric is not compact. 
\end{proof}
\end{lemma}

\begin{remark}[Conditions on weak lower semicontinuity]
\label{rmk_weakly_conti}

    We give several examples of functionals satisfying weak lower semicontinuity.
    \begin{itemize}[leftmargin=*]

        \item If $f:\Theta\rightarrow(-\infty,+\infty]$ is lower semicontinuous and its negative part has a 2-growth, then $F=\int_{\Theta} fd\rho$ is lower semicontinuous with respect to $\mathcal{W}_2$ topology, see Example 9.3.1 in \cite{ambrosio2005gradient}. If $f:\Theta\rightarrow(-\infty,+\infty]$ is lower semicontinuous and bounded from below, then $F=\int_{\Theta} fd\rho$ is weakly lower semicontinuous, see Example 9.3.1 in \cite{ambrosio2005gradient}. 
        \item The KL divergence $D_{KL}(\rho\Vert\tilde\rho)$ is jointly weakly lower semicontinuous, see Lemma 9.4.3 in \cite{ambrosio2005gradient}.
    \end{itemize}
\end{remark}

\begin{remark}[Conditions on (1) in Assumption \ref{assum_1}] 
\label{rmk_admit_minimizer}
% (a) If $F:\mathcal{P}_2(\Theta):\rightarrow (-\infty,+\infty]$ is proper, lower semicontinuous (with respect to $\mathcal{W}_2$ topology), and bounded from below, then $F$ is weak lower semicontinuous.

If $F:\mathcal{P}_2(\Theta)\rightarrow (-\infty,+\infty]$ is proper, lower 
semicontinuous (with respect to $\mathcal{W}_2$ topology), $\mu$-convex ($\mu\in\mathbb{R}$) along generalized geodesics, and bounded from below, then the proximal algorithm (\ref{Eqn_jko}) admits a minimizer for any $\xi>0$. For details, we refer to Section~10.3 in~\cite{ambrosio2005gradient}, where the condition is even weaker, being defined only on $\mathcal{W}_2$-bounded sets.

% If $f:\Theta\rightarrow(-\infty,+\infty]$ is lower semicontinuous, and $\mu$-convex ($\mu\in\mathbb{R}$) in Euclidean sense, then $F=\int_{\Theta} fd\rho$ is weak lower semicontinuous.
\end{remark}

% \begin{proof}[Proof of Lemma~\ref{lemma_weakly_conti}]
%     See Section 10.3 in \cite{ambrosio2005gradient}.
% \end{proof}
\begin{lemma}[Satisfaction of Assumption \ref{assum_2} on $\mathcal{P}_2(\Theta)$ with compact set $\Theta$]
 \label{lemma_compact} Under Assumption \ref{assum_1}, if $\Theta$ is compact, then Assumption \ref{assum_2} holds.
 \begin{proof} By Proposition 7.17 in~\cite{santambrogio2015optimal}, when $\Theta$ is compact, for fixed $ \rho\in \mathcal{P}_2(\Theta)$, the first variation of $\mathcal{W}_2(\cdot,\rho)$ is well-defined for all $\tilde\rho\in\mathcal{P}_2^a(\Theta)$. Similar to Proposition 8.7 in \cite{santambrogio2015optimal}, by standard calculus of variation followed by gradient operation, we have
 \begin{equation}
     \nabla\dfrac{\delta F}{\delta \rho}(\rho_\xi)=\dfrac{T_{\rho_\xi}^{\rho}-\mathbf{id}}{\xi}.
 \end{equation}
 Please refer to Lemma B.1 in \cite{yao2024wasserstein} for similar arguments. Therefore, Assumption \ref{assum_2} holds. 
 \end{proof}
\end{lemma}
% \begin{roundedbox}
%     \begin{assumption}[Smoothness assumption, {[Assumption 1, \cite{chizat2022meanfield}]}]
%         \label{assum_3}

% \end{assumption}
% \end{roundedbox}

\begin{lemma} 
\label{integral_lemma}Assume $u:[0,\xi]\rightarrow (0,+\infty)$ is a decreasing function, $\partial_t u(t)\leq -C(t) u(t)$ almost everywhere for $t\in[0,\xi]$, $C(t)>0$ for every $t\in[0,\xi]$, then 
\begin{equation}
    u(\xi)\leq u(0)\exp\left(\int_0^{\xi}-C(t)d t \right).
\end{equation}

\begin{remark}
    Lemma~\ref{integral_lemma} extends the classical Gronwall lemma, which requires everywhere differentiability, to the case with almost everywhere differentiability and monotonicity.  
\end{remark}

\begin{proof}
    By the monotonicity of $u(t)$, we construct a function $g(t)=\ln(u(t))$. It is a decreasing function and $\dfrac{dg(t)}{dt}=\dfrac{\partial_t u(t)}{u(t)}\leq -C(t)$ almost everywhere on $[0,\xi]$.
 By properties of the Lebesgue integral, we have
\begin{equation}
    \int_0^{\xi} \dfrac{dg(t)}{dt} dt \leq \int_0^{\xi}-C(t)dt.
\end{equation} Since $g$ is decreasing, by Proposition 6.6 in \cite{komornik2016lectures},
{\begin{equation}
    g(\xi)-g(0)\leq\int_0^{\xi} \dfrac{dg(t)}{dt} dt.
\end{equation}}
Thus,
\begin{equation}
    u(\xi)\leq u(0) \exp\left(\int_0^{\xi}-C(t)dt \right).
\end{equation}
\end{proof}
% \begin{proof}
%     By the monotonicity of $u(\xi)$, we have the integral form of Gronwall Formula,
% \begin{equation}
%     u(\xi)\leq u(0)+
% \end{equation}
% \end{proof}

\end{lemma}
% \begin{lemma} 
% \label{lemma_otmap}Assume $T_{\overline{\rho}_n}^{\overline{\rho}_{n+1}}$ is $C^2$ diffeomorphism. If  $\overline{\rho}_n\in\mathcal{C}^1(\Theta)$, then $\overline{\rho}_{n+1}\in \mathcal{C}^1(\Theta)$.

% \begin{proof}
% The change variable formula of probability density is, 
% \begin{equation}
%     T_{\#}{\rho}(\theta)=\rho(T^{-1}(\theta))\cdot |\det D(T^{-1})(\theta)|,
% \end{equation}
% where $D(T^{-1})$ is the Jacobian matrix of $T^{-1}$.
%     Since $T$ is a $\mathcal{C}^2$ diffeomorphism, then $T^{-1}$ is $\mathcal{C}^2$ mapping. Thus, $\rho\circ (T^{-1})$ is $\mathcal{C}^1(\Theta)$. Since $T^{-1}$ is diffeomorphism, then $D(T^{-1})$ is not singular and $|\det D(T^{-1})|$ is $C^1(\Theta)$. Then $T_{\#}\rho$ is $\mathcal{C}^1(\Theta)$.
% \end{proof}
% \end{lemma}
\begin{remark}
\label{ref:realnvp} The coupling layer in Real NVP \citep{dinh2016density} has the following structure,
\begin{equation}
T:\left\{
\begin{aligned}
&\alpha_{1:\tilde{d}}\quad       = \theta_{1:\tilde{d}}, \\
&\alpha_{\tilde{d}+1:d}     = \theta_{\tilde{d}+1:d} * \exp(s(\theta_{1:\tilde{d}})) + t(\theta_{1:\tilde{d}}),
\end{aligned}
\right .
\end{equation}
 where $*$ refers to the pointwise product. 
\end{remark}
% \begin{lemma}
% \label{lemma_C1}
%  Under Assumption \ref{assum_4}, if $\overline{\rho}_n\in\mathcal{C}^1(\Theta)$ and $\overline\rho_{n+1}$ is the exact solution of the Wasserstein proximal algorithm (\ref{Eqn_jko}), then $\beta_{n+1}=0$. 
%  \begin{proof}
%      Since $\overline{\rho}_{n+1}$ is the exact solution of (\ref{Eqn_jko}), then $\overline{\rho}_{n+1}\in D(|\partial F|)$ and $\frac{T_{\overline{\rho}_{n+1}}^{\overline\rho_n}-\mathbf{id}}{\xi}\in \partial F(\overline{\rho}_{n+1})$ by Lemma 10.1.2 in \cite{ambrosio2005gradient}. As $\overline{\rho}_{n+1}\in \mathcal{C}^1(\Theta)$, then $\nabla \dfrac{\delta F}{\delta \rho}(\overline{\rho}_{n+1})$ is the unique subdifferential by Lemma 10.4.1 in \cite{ambrosio2005gradient}. Therefore, $${T_{\overline{\rho}_{n+1}}^{\overline\rho_n}-\mathbf{id}}-\xi\nabla \frac{\delta F}{\delta \rho}(\overline{\rho}_{n+1})=0.$$
%  \end{proof}
% \end{lemma}

\section{Proofs}
\label{app:proofs}

%The proof of Lemma \ref{Lemma_important} essentially follows from Proposition 3.1 and Proposition 3.3 in \citep{Ambrosio_2013}, which are summarized together in Lemma \ref{lem:ambrosio_tech} for completeness. 

\begin{proof}[{\bf Proof of Lemma \ref{Lemma_important}}]
    By Assumption \ref{assum_1}, the proximal problem admits a minimizer. Thus we can define
    $$z_\xi\in \arg\min_{z'\in\mathcal P_2(\Theta)}
\left\{
F(z')+\frac{1}{2\xi}\mathcal W_2^2(z',z)
\right\}.$$
Then Lemma~\ref{lem:ambrosio_tech} implies that $D_{+}(z,\xi)=D_{-}(z,\xi)=d(z,z_\xi)$ in $(0,t^*(z))$ except for at most countable exceptions, from which we can conclude the desired Lemma \ref{Lemma_important}.
\end{proof}

\begin{proof}[Proof of Theorem~\ref{main_thm}] It suffices to prove one-step contraction. We begin with
\begin{align*}
    \partial_{\xi}(u(\rho,\xi)-F^*)&= -\dfrac{\mu}{2\xi(1+\mu \xi)}\mathcal{W}_2^2(\rho_\xi,\rho)-\dfrac{1}{2(1+\mu\xi)\xi^2}\mathcal{W}_2^2(\rho_\xi,\rho) & \text{(by Lemma~\ref{Lemma_important})}\\
    &=-\dfrac{\mu}{2\xi(1+\mu \xi)}\mathcal{W}_2^2(\rho_\xi,\rho)-\dfrac{1}{2(1+\mu \xi)\xi^2}\int_{\Theta} \left\Vert T_{\rho_{\xi}}^{\rho}-\mathbf{id}\right\Vert^2 d\rho_{\xi} &\text{(by Brenier's theorem)}\\
    &\leq  -\dfrac{\mu}{2\xi(1+\mu \xi)}\mathcal{W}_2^2(\rho_\xi,\rho)-\dfrac{1}{2(1+\mu \xi)}\int_{\Theta} \left\Vert \nabla \dfrac{\delta F}{\delta \rho}(\rho_\xi)\right\Vert^2 d\rho_\xi & \text{(by Assumption \ref{assum_2})}\\
    &\leq -\dfrac{\mu}{2\xi(1+\mu \xi)}\mathcal{W}_2^2(\rho_\xi,\rho)-\dfrac{\mu}{(1+\mu \xi)}(F(\rho_\xi)-F^{*}) &\text{(by PL (\ref{Eqn_wass_PL}))}\\
    &=-\dfrac{\mu}{1+\mu\xi}(u(\rho,\xi)-F^*).& \text{(by Def \ref{Def_hopf})}
\end{align*} 
 Using Lemma \ref{integral_lemma} to deal with the technical issue of almost everywhere differentiability, we have
\begin{align*}
    u(\rho,\xi)-F^* &\leq (u(\rho,0)-F^*)\exp\left(\int_0^{\xi} -\frac{\mu}{1+\mu t}dt\right)=(F(\rho)-F^*)\dfrac{1}{1+\mu\xi}.
\end{align*}
Invoking the PL inequality (\ref{Eqn_wass_PL}) once again, we obtain that
\begin{align*}
    (F(\rho)-F^*)\dfrac{1}{1+\mu\xi}&\geq u(\rho,\xi)-F^*=F(\rho_\xi)-F^*+\dfrac{1}{2\xi}\mathcal{W}_2^2( \rho_\xi,\rho)\\&\geq F(\rho_\xi)-F^*+\dfrac{\xi}{2}\int_{\Theta} \left\Vert \nabla \dfrac{\delta F}{\delta \rho}(\rho_\xi)\right\Vert^2 d\rho_\xi\geq (1+\mu\xi)(F(\rho_\xi)-F^*).
\end{align*}
\end{proof}

\begin{definition}[Uniform log-Sobolev inequality]
\label{def:unifLSI}
    There is a constant $\mu > 0$ such that for any $\rho \in \mathcal{P}_2(\Theta)$, its Gibbs proximal distribution $q_\rho$ defined as, 
   \begin{equation}
    \label{eqn:proxGibbs}
    q_\rho(\theta) \propto \exp\left( -{1 \over \tau} {\delta R\over \delta \rho}(\rho) (\theta) \right)
\end{equation}
satisfies the log-Sobolev inequality (\ref{eqn:LSI}) with the constant $\mu$. 

\end{definition}
\begin{proof}[{\bf Proof of Corollary \ref{Cor_mfld}}] We divide our proof into four parts. \\
(1) We prove the satisfaction of Assumption \ref{assum_1} and the geodesic semiconvexity of $F_{\tau}$.\\
(2) We prove the PL inequality.\\
(3) We prove the satisfaction of Assumption \ref{assum_2} by showing that $\nabla\frac{\delta F_\tau}{\delta \rho}(\rho_\xi)=\partial^{\circ}F_\tau(\rho_\xi)$.\\
(4) With the previous three parts, we can get the linear convergence of the function value. The last part is devoted to obtaining a convergence rate of $\mathcal{W}_2$ distance using some structure of MFLD.

\textbf{Our proof is as follows.}\\
(1) The weak lower semicontinuity and lower boundedness of $F_{\tau}$ are verified in Section 5.1 in \cite{chizat2022meanfield}, and thus a minimizer is admitted by Lemma \ref{Lemma_exist}. The geodesic semiconvexity follows from Lemma A.2 in \cite{chizat2022meanfield} and the proof relies on (\ref{eqn_assum_3}) below.

(2) Now we prove the PL inequality. The training risk $R:\mathcal{P}_2(\Theta)\rightarrow (-\infty,+\infty]$ has linear convexity if the loss function $l$ is convex (in the Euclidean sense). By Proposition 5.1 in~\citep{chizat2022meanfield}, the $L^2$-regularized training risk $R$ in (\ref{twonn}) satisfies the uniform LSI assumption (Assumption 3 in \cite{chizat2022meanfield}). Next, we shall show that these assumptions imply the relaxed PL-inequality defined in (\ref{Eqn_wass_PL}). By the entropy sandwich bound Lemma 3.4 in~\citep{chizat2022meanfield} (which relies on linear convexity of $F_\tau$),
        \begin{equation}
         \tau \KL(\rho \| q_\rho) \geq F_\tau(\rho)-F_\tau(\rho^*).
        \end{equation}
        Therefore,
        \begin{align*}
           \int_{\Theta} \left\Vert\nabla \frac{\delta F_\tau}{\delta \rho}(\rho)\right\Vert^2 d\rho&=\int_{\Theta} \left\Vert \nabla\frac{\delta{R}}{\delta \rho}(\rho)+\tau\nabla \log(\rho)\right\Vert^2 d\rho\\
           &=\tau^2 J_{q_\rho}(\rho) \geq 2\mu_\tau\tau^2 \KL(\rho \| q_\rho) \geq 2\mu_\tau \tau(F_\tau(\rho)-F_\tau(\rho^*)).
        \end{align*} 
Thus, the functional $F_\tau$ satisfies the Wasserstein PL-inequality with parameter $\tau \mu_\tau$.

(3) Under Assumption \ref{assum_1}, $(T_{\rho_\xi}^{\rho}-\mathbf{id}) / \xi$ is a Fr\'echet subdifferential at $\rho_\xi\in\mathcal{P}_2^a(\Theta)$, and $\rho_{\xi}\in D(|\partial F_\tau|)$ by Lemma 10.1.2 in \cite{ambrosio2005gradient}. To prove that Assumption \ref{assum_2} holds, we only need to prove that, 
\begin{center}
 \textit{If $\rho\in D(|\partial F_\tau|)\cap \mathcal{P}_2^a(\Theta)$, then $\nabla \frac{\delta F_\tau}{\delta \rho}(\rho)= \partial^{\circ} F_\tau(\rho)$. }
 \end{center}

Our proof for (3) below highly relies on the proof of  Theorem 10.4.13 in \cite{ambrosio2005gradient}.\\
\textbf{Step 1.} We first need to derive some conditions similar to equations (10.4.58) and (10.4.59) in \cite{ambrosio2005gradient}. Under the assumptions of Corollary \ref{Cor_mfld}, $R$ satisfies the following smoothness condition with $\tilde{L}>0$ by Proposition 5.1 in \cite{chizat2022meanfield}. Thus, $\forall \theta,\tilde\theta\in\Theta$ and $\forall \rho,\tilde\rho\in\mathcal{P}_2(\Theta)$, we have 
\begin{equation}
\label{eqn_assum_3}
      \left\Vert \nabla \dfrac{\delta{R}}{\delta \rho}(\rho)(\theta)- \nabla\dfrac{\delta{R}}{\delta \rho}(\tilde\rho)(\tilde\theta)\right\Vert \leq \tilde{L}(\Vert \theta-\tilde\theta\Vert_2+\mathcal{W}_2(\rho,\tilde\rho)).
\end{equation}
By Lemma A.2 in \cite{chizat2022meanfield}, relying on (\ref{eqn_assum_3}), choosing $\mathbf{r}=\mathbf{id}+\mathbf{t}$ with $\mathbf{t}\in C_c^{\infty}(\Theta)$,
\begin{equation} 
\label{Eqn_lemma_A2}
    \lim_{t\rightarrow 0}\dfrac{R((\mathbf{id}+t\mathbf{t})_{\#}\rho)-R(\rho)}{t}=\int_{\Theta} \left\langle \nabla \dfrac{\delta R}{\delta \rho}(\rho), \mathbf{r}-\mathbf{id} \right\rangle d\rho,
\end{equation}
which is the key condition similar to equations (10.4.58) and (10.4.59) in \cite{ambrosio2005gradient} that we want to obtain. 

Furthermore, a by-product is that (\ref{Eqn_lemma_A2}) suggests that $-\nabla\frac{\delta R}{\delta \rho}(\rho)$ is a subdifferential of $(-R)$ and $\left\Vert\nabla\frac{\delta R}{\delta \rho}(\rho)\right\Vert_{L_2(\rho)}=\left|\partial (-R)\right|(\rho)$ for all $\rho\in \mathcal{P}_2^a(\Theta)$, see Lemma A.2 in \cite{chizat2022meanfield}. It is not hard to verify that $\left\Vert\nabla\frac{\delta R}{\delta \rho}(\rho)\right\Vert_{L_2(\rho)}$ is finite at any $\rho\in\mathcal{P}_2(\Theta)$ because (\ref{eqn_assum_3}) ensures 2-growth of $\left\Vert\nabla \frac{\delta R}{\delta \rho}(\rho)(\cdot)\right\Vert^2$.

In \textbf{Step 2}, we conduct a proof similar to Theorem 10.4.13 in \cite{ambrosio2005gradient}.

\textbf{Step 2.} By Lemma 10.4.4 in \cite{ambrosio2005gradient}, (\ref{Eqn_lemma_A2}), 
%the inequality $\mathcal{W}_2((\mathbf{i}+t\mathbf{t})_{\#}\rho,\rho)\leq \Vert\mathbf{t}\Vert_{L_2(\rho)}$ ($t\in(0,1)$), 
and the fact that $\rho\in \mathcal{P}_2^a(\Theta)\cap D(|\partial F_{\tau}|)$, for $\mathbf{t}\in C_{c}^{\infty}(\Theta)$
\begin{equation}
    -\int_{\Theta} \tau\rho\nabla\cdot \mathbf{t}d\theta+\int_{\Theta} \left\langle \nabla \dfrac{\delta R}{\delta \rho}(\rho), \mathbf{t} \right\rangle d\rho\geq -\vert\partial F_\tau(\rho)\vert \left\Vert \mathbf{t}\right\Vert_{L_2(\rho)}.
\end{equation}
 By (\ref{eqn_assum_3}), $  \forall \theta,\tilde\theta\in\Theta, \forall{\rho}\in\mathcal{P}_2(\Theta)$, 
\begin{equation}
\label{eqn_bound}
  \quad \left\Vert \nabla \dfrac{\delta{R}}{\delta \rho}(\rho)(\theta)- \nabla\dfrac{\delta{R}}{\delta \rho}(\rho)(\tilde\theta)\right\Vert \leq \tilde L\Vert \theta-\tilde\theta\Vert_2.
\end{equation}
Therefore, $\nabla \dfrac{\delta{R}}{\delta \rho}(\rho)(\cdot)$ is Lipschitz and locally bounded. Therefore,  following the
 same argument of Theorem 10.4.13 in \cite{ambrosio2005gradient}, we obtain that 
 \begin{equation}
     \nabla \dfrac{\delta F_\tau}{\delta \rho}(\rho)=\tau\dfrac{\nabla \rho}{\rho}+\nabla \dfrac{\delta R}{\delta \rho}(\rho)\in L_2(\rho) \text{ and } {\left\Vert\nabla\dfrac{\delta F_\tau}{\delta \rho}(\rho)\right\Vert_{L_{2}(\rho)}\leq |\partial F_\tau(\rho)|}.
 \end{equation} where we define  $\frac{\nabla \rho}{\rho}=0$ if $\rho(\theta)=0$. Furthermore, it is straightforward to prove $\nabla\dfrac{\delta F_\tau}{\delta \rho}(\rho)\in\partial F_\tau(\rho)$.\footnote{Here the proof is slightly different from Theorem 10.4.13 in \cite{ambrosio2005gradient} for simplicity,  as we already have the finite slope property of $R$.}
 Since $\rho\in D(|\partial F_\tau|)$ and $\rho\in D(|\partial (-R)|)$, we have that $|\partial H_\tau|(\rho)\leq |\partial F_\tau|(\rho)+|\partial (-R)|(\rho)$ where $H_\tau(\rho)=\int_{\Theta} \tau\rho \log\rho$. Therefore, we have $\rho\in D(|\partial H_\tau|)$. Thus, $\tau\frac{\nabla{\rho}}{\rho}\in \partial H_\tau(\rho)$ in view of Theorem 10.4.6 in \cite{ambrosio2005gradient}. With $\nabla\dfrac{\delta R}{\delta \rho}(\rho)\in \partial R(\rho)$, we have proved that {$\nabla\dfrac{\delta F_\tau}{\delta \rho}(\rho)\in\partial F_\tau(\rho)$.} Thus, $\nabla\dfrac{\delta F_\tau}{\delta \rho}(\rho)=\partial^{\circ} F_\tau(\rho)$.
 
(4) With the previous three parts, we can get the first inequality in Corollary \ref{Cor_mfld} with Theorem \ref{main_thm}. Next, we prove the second inequality in Corollary \ref{Cor_mfld}. Using Lemma 3.4 in \citep{chizat2022meanfield} once again, we get
$$\tau \KL(\rho \| \rho^*)\leq F_{\tau}(\rho)-F_\tau(\rho^*)\leq \tau \KL(\rho \| q_\rho).$$
Since $\rho^*$ also satisfies $\mu_\tau$-LSI, by Talagrand's inequality, we have
\begin{equation}
    \mathcal{W}_2^2(\rho^*,\rho_n)\leq \dfrac{2}{\mu_\tau} \KL(\rho_n \| \rho^*)\leq \dfrac{2}{\mu_\tau \tau}(F_\tau(\rho_n)-F_\tau(\rho^*))\leq  \dfrac{2}{\mu_\tau \tau}(F_\tau(\rho_0)-F_\tau(\rho^*))\left(\dfrac{1}{1+\mu_\tau\tau\xi}\right)^{2n}.
\end{equation}
\end{proof}

\begin{proof}[{\bf Proof of Corollary \ref{cor_ld}}] Similar to proof of Corollary \ref{Cor_mfld}, we divide the proof into four parts. 

(1) $F$ is weakly lower semicontinuous by Remark \ref{rmk_weakly_conti} and is bounded from below because of the LSI condition. By Remark \ref{rmk_admit_minimizer}, $prox_{F,\xi}$ admits one minimizer in $\mathcal{P}_2(\Theta).$ In addition, $D(F)\subset \mathcal{P}_2^a(\Theta)$. Thus, Assumption \ref{assum_1} is satisfied.

(2) The PL inequality directly follows from the LSI condition. 

(3) We prove that Assumption \ref{assum_2} is satisfied by showing  $\nabla\frac{\delta F}{\delta \rho}(\rho_\xi)=\partial^{\circ}F(\rho_\xi)$. Theorem 10.4.13 in \cite{ambrosio2005gradient} already incorporates KL divergence as a special case: Assume $f$ is semiconvex and lower semicontinuous, and $\{\theta|f(\theta)<\infty\}$ is not empty. Then for $F(\rho)=\int_{\Theta} f d\rho+\int_{\Theta} \rho \log\rho$ (here we assume the interaction energy to be 0), if $\rho\in D(|\partial F|)\cap \mathcal{P}_2^a(\Theta)$, then $\nabla \frac{\delta F}{\delta \rho}(\rho)=\nabla f+\nabla \log(\rho)=\partial^{\circ}F(\rho)$, where we assume $\nabla \log(\rho)(\theta)=0$ if $\rho(\theta)=0$. With Assumption \ref{assum_1} to guarantee $\rho_\xi\in\mathcal{P}_2^a(\Theta)$ and Lemma 10.1.2 in \cite{ambrosio2005gradient} to guarantee that $\rho_\xi\in D(|\partial F|)$, Assumption \ref{assum_2} holds.

(4) The convergence rate on $\mathcal{W}_2$ distance follows Talagrand's inequality. 
\end{proof}
\begin{proof}[{\bf Proof of Example \ref{Cor_specialexample}}]
(1) We first verify Assumption \ref{assum_1}. By definition,
\(D(F)=\mathcal P_2^a(\mathbb R)\) since $\int |\theta|d\rho(\theta)\leq \left (\int \theta^2 d\rho(\theta)\right )^{1/2}< \infty$. Fix \(\rho\in D(F)\), \(\xi>0\), and let $
    x=m(\rho)$ and
   $$ H(m)=m+c\sin m. $$
For any \(\tilde\rho\in D(F)\), denote \(\tilde x=m(\tilde\rho)\). By Jensen's inequality,
$
    W_2^2(\tilde\rho,\rho)
    \ge
    |m(\tilde\rho)-m(\rho)|^2
    =
    |\tilde x-x|^2.
$
Conversely, for every \(\tilde x\in\mathbb R\), define
\[
    \tilde\rho_{\tilde x}=(\mathbf{id}+\tilde x-x)_\#\rho.
\]
Since \(\rho\in\mathcal P_2^a(\mathbb R)\), we have
\(\tilde\rho_{\tilde x}\in\mathcal P_2^a(\mathbb R)\). Moreover,
$
    m(\tilde\rho_{\tilde x})=\tilde x,
    W_2^2(\tilde\rho_{\tilde x},\rho)=|\tilde x-x|^2 $.
Therefore, 
\[
    \inf_{\tilde\rho\in\mathcal P_2(\mathbb R)}
    \left\{
        F(\tilde\rho)+\frac{1}{2\xi}W_2^2(\tilde\rho,\rho)
    \right\}=   \inf_{\tilde x\in\mathbb R}
    \left\{
        H(\tilde x)^2+\frac{1}{2\xi}|\tilde x-x|^2
    \right\}.
\]
The RHS is continuous, and coercive, i.e.,$H(\tilde x)^2+1/{2\xi}|\tilde{x}-x|^2\rightarrow \infty$ when $|\tilde x|\rightarrow \infty$, hence admits a minimizer
$x_\xi$. Therefore, $\rho_\xi=(\mathrm{id}+x_\xi-x)_\#\rho$
is a proximal minimizer and belongs to \(D(F)\). Thus Assumption
\ref{assum_1} is satisfied.

\noindent
(2) We next verify the Wasserstein PL inequality. On \(D(F)\), we have
\[
    F(\rho)=H(m(\rho))^2.
\]
Thus,
\[
    \nabla \frac{\delta F}{\delta\rho}(\rho)(\theta)
    =
    2H(m(\rho))H'(m(\rho)).
\]
Since $ H'(m)=1+c\cos m\ge 1-c>0 $,
we obtain
\[
    \int_{\mathbb R}
    \left\Vert
        \nabla \frac{\delta F}{\delta\rho}(\rho)(\theta)
    \right\Vert^2
    d\rho(\theta)
    =
    4H(m(\rho))^2H'(m(\rho))^2
    \ge
    4(1-c)^2F(\rho).
\]
Moreover, \(F^\ast=0\), for example by taking any
\(\rho\in\mathcal P_2^a(\mathbb R)\) with \(m(\rho)=0\). Therefore,
\[
    \int_{\mathbb R}
    \left\Vert
        \nabla \frac{\delta F}{\delta\rho}(\rho)(\theta)
    \right\Vert^2
    d\rho(\theta)
    \ge
    2\mu \big(F(\rho)-F^\ast\big),
    \qquad
    \mu=2(1-c)^2.
\]
Thus the Wasserstein PL inequality is satisfied.

\noindent
(3) We now verify Assumption \ref{assum_2}. Since $\rho$ is the proximal point, then $\rho\in D(|\partial F|)$ and $\partial F(\rho)$ is not empty. Let \(v= \partial^\circ F(\rho)\),  for any
\(\mathbf{t}\in C_c^\infty(\mathbb R)\), consider the perturbation
\[
    \rho_t = (\mathbf{id}+t\mathbf{t})_\#\rho .
\]
For $t>0$ sufficiently small, \(\mathbf{id}+t\mathbf{t}\)
is an increasing rearrangement, and hence it is the optimal transport map from \(\rho\) to
\(\rho_t\). By the definition of the subdifferential,
\[
    F(\rho_t)-F(\rho)
    \ge
    t\int_{\mathbb R} v(\theta)\mathbf{t}(\theta)\,d\rho(\theta)
    +o(\mathcal W_2(\rho,\rho_t))=t\int_{\mathbb R} v(\theta)\mathbf{t}(\theta)\,d\rho(\theta)
    +o(t\Vert \mathbf{t}\Vert_{L_2(\rho)}).
\]
On the other hand,
\[
    m(\rho_t)
    =
    \int_{\mathbb R}(\theta+t\mathbf{t}(\theta))\,d\rho(\theta)
    =
    m(\rho)+t\int_{\mathbb R}\mathbf{t}(\theta)\,d\rho(\theta).
\]
Thus
\[
    F(\rho_t)-F(\rho)
    =
    2H(m(\rho))H'(m(\rho))
    t\int_{\mathbb R}\mathbf{t}(\theta)\,d\rho(\theta)
    +o(t \Vert \mathbf{t}\Vert_{L_2(\rho)}).
\]
Dividing by \(t>0\) and letting \(t\downarrow0\), and then
repeating the argument with $-\mathbf{t}$, gives
\[
    \int_{\mathbb R}v(\theta)\mathbf{t}(\theta)\,d\rho(\theta)
    =
    \int_{\mathbb R}
    2H(m(\rho))H'(m(\rho))\mathbf{t}(\theta)\,d\rho(\theta).
\]
Since \(C_c^\infty(\mathbb R)\) is dense in \(L^2(\rho)\), we obtain
\[
    v(\theta)
    =
    2H(m(\rho))H'(m(\rho))
    \qquad \rho\text{-a.e.}
\]
Thus Assumption \ref{assum_2} holds.

\noindent
(4) Finally, we show that \(F\) is neither linearly convex nor geodesically
convex when $c=1/2$. Take two measures $\rho_0=\mathcal{N}(\pi/2,1), \rho_1=\mathcal{N}(3\pi/2,1)$. And thus we have $m(\rho_0) = \frac{\pi}{2}, \quad m(\rho_1) = \frac{3\pi}{2}.$ Then for the linear interpolation $\rho_t = (1-t)\rho_0 + t\rho_1$, we have $m(\rho_t) = \frac{\pi}{2} + \pi t.$ Let $G(m) = H^2(m)= (m + \frac{1}{2} \sin m)^2.$ Then $$ F(\rho_t) = G\left(\frac{\pi}{2} + \pi t\right). $$

A direct computation shows $G''(m)=2(1+\frac{1}{2} \cos m)^2-\sin m (m+\frac{1}{2} \sin m)$ and $G''\left(\frac{\pi}{2}\right) < 0,$ so $t \mapsto F(\rho_t)$ is not convex. Hence $F$ is not linearly convex. Since the optimal transport map between $\rho_0$ and $\rho_1$ is nothing but translation, Wasserstein geodesics between $\rho_0$ and $\rho_1$ preserve the mean interpolation: $m(\tilde\rho_t) = (1-t)m(\rho_0) + t m(\rho_1)=\frac{\pi}{2} + \pi t$, where $\tilde \rho_t$ is the push-forward measure of $\rho_0$ by $(1-t)\mathbf{id}+tT$ with optimal map $T(\theta)=\theta+\pi$. Thus along the $W_2$ geodesic, $F(\tilde\rho_t) = G\left(\frac{\pi}{2} + \pi t\right),$ is also not geodesically convex.

\end{proof}
\begin{proof}[{\bf Proof of Theorem \ref{Cor_sharp}}]In this proof, {we will not invoke Assumption \ref{assum_2}.}\\
\textbf{Step 1.} We want to prove that $\mu$-geodesic convex implies, for any fixed $\rho$ and for any $\xi$, 
\begin{equation}
\label{Eqn_strong_invoke}
    F(\rho^*)\geq F(\rho_\xi)-\dfrac{1}{2\mu\xi^2}\mathcal{W}_2^2(\rho_\xi,\rho).
\end{equation}
Note $\dfrac{T_{\rho_\xi}^{\rho}-\mathbf{id}}{\xi}$ is a Fr\'echet subdifferential at $\rho_\xi$,
\begin{align*}
 F(\tilde{\rho})&\geq F(\rho_\xi)+\int_{\Theta}\left\langle \dfrac{T_{\rho_\xi}^{\rho}(\theta)-\theta}{\xi} ,T_{\rho\xi}^{\tilde{\rho}}(\theta)-\theta \right\rangle d\rho_\xi(\theta)+\dfrac{\mu}{2}{\mathcal{W}_2^2(\rho_\xi,\tilde{\rho})}\\ 
 &=F(\rho_\xi)+\int_{\Theta}\left(\left\langle \dfrac{T_{\rho_\xi}^{\rho}(\theta)-(\theta)}{\xi} ,T_{\rho_\xi}^{\tilde{\rho}}(\theta)-\theta \right\rangle+\dfrac{\mu}{2}\left\Vert T_{\rho_\xi}^{\tilde{\rho}}(\theta)-\theta\right\Vert^2 \right)d\rho_\xi(\theta).
\end{align*}
%  Then, we minimize both sides of (\ref{eqn:mu_conxity}) with respect to $\tilde{\rho}\in\mathcal{P}_2^a(\Theta)$. Clearly, $\rho^*$ minimizes the left side.  For the integral term on the right-hand side, we define $T_{\rho_\xi}^{\tilde{\rho}}$ in the following,
%  \begin{equation}
%      T_{\rho_\xi}^{\tilde{\rho}}(\theta)-\theta=-\dfrac{1}{\mu\xi}\left(T_{\rho_\xi}^{\rho}(\theta)-\theta\right), \quad \rho_\xi\text{-}a.e.
% \label{eqn_sharp_reb}
%  \end{equation}
% then it minimizes the integral as it minimizes the term inside the integral almost everywhere. Therefore,
% \begin{equation}
%      F(\rho^*)\geq F(\rho_\xi)-\dfrac{1}{2\mu\xi^2}\int_{\Theta}\left\Vert T_{\rho_\xi}^{\rho}(\theta)-\theta\right\Vert^2 d\rho_{\xi}(\theta)=F(\rho_\xi)-\dfrac{1}{2\mu\xi^2}W_2^2(\rho_\xi,\rho).
% \end{equation}
{Then, taking $\tilde{\rho}=\rho^*$ in (\ref{eqn:mu_conxity}), we obtain
\begin{equation}
\begin{aligned}
F(\rho^*)\geq F(\rho_\xi)
&+
\int_{\Theta}
\left\langle
\frac{T_{\rho_\xi}^{\rho}(\theta)-\theta}{\xi},
T_{\rho_\xi}^{\rho^*}(\theta)-\theta
\right\rangle d\rho_\xi(\theta)+
\frac{\mu}{2}W_2^2(\rho_\xi,\rho^*).
\end{aligned}
\end{equation}
By Cauchy--Schwarz, we have
$$
\int_{\Theta}
\left\langle
\frac{T_{\rho_\xi}^{\rho}(\theta)-\theta}{\xi},
T_{\rho_\xi}^{\rho^*}(\theta)-\theta
\right\rangle d\rho_\xi(\theta)
\geq
-\frac{1}{\xi}W_2(\rho_\xi,\rho)W_2(\rho_\xi,\rho^*).
$$
Therefore,
\begin{equation}
F(\rho^*)\geq F(\rho_\xi)
-\frac{1}{\xi}W_2(\rho_\xi,\rho)W_2(\rho_\xi,\rho^*)
+
\frac{\mu}{2}W_2^2(\rho_\xi,\rho^*)\geq F(\rho_\xi)
     -\dfrac{1}{2\mu\xi^2}W_2^2(\rho_\xi,\rho).
\end{equation}}

\textbf{Step 2.} \begin{align*}
    \partial_{\xi}(u(\rho,\xi)-F^*)&= -\dfrac{\mu}{2\xi(1+\mu \xi)}\mathcal{W}_2^2(\rho_\xi,\rho)-\dfrac{1}{2(1+\mu\xi)\xi^2}\mathcal{W}_2^2(\rho_\xi,\rho) & \text{(by (\ref{Eqn_important_lemma}))}\\
    &\leq -\dfrac{\mu}{2\xi(1+\mu \xi)}\mathcal{W}_2^2(\rho_\xi,\rho)-\dfrac{\mu}{(1+\mu \xi)}(F(\rho_\xi)-F^{*}) &\text{(by (\ref{Eqn_strong_invoke}) in Step 1)}\\
    &=-\dfrac{\mu}{1+\mu\xi}(u(\rho,\xi)-F^*).& \text{(by Def \ref{Def_hopf})}
\end{align*}
Using Lemma \ref{integral_lemma} to deal with the technical issue of almost everywhere differentiability, we have
\begin{align*}
    u(\rho,\xi)-F^* &\leq (u(\rho,0)-F^*)\exp\left(\int_0^{\xi} -\frac{\mu}{1+\mu t}dt\right)=(F(\rho)-F^*)\dfrac{1}{1+\mu\xi}.
\end{align*}
Invoking (\ref{Eqn_strong_invoke}) once again, we obtain that
\begin{align*}
    (F(\rho)-F^*)\dfrac{1}{1+\mu\xi}\geq u(\rho,\xi)-F^*=F(\rho_\xi)-F^*+\dfrac{1}{2\xi}\mathcal{W}_2^2( \rho_\xi,\rho)\geq (1+\mu\xi)(F(\rho_\xi)-F^*).
\end{align*}
\textbf{Step 3.} We derive a bound for $\mathcal{W}_2(\rho^*,\rho_n)$ here. Since for any $\rho\in \mathcal{D}(F)$, $$F(\rho)-F(\rho^*)\geq \int_{\Theta}\langle {\bf{0}}, T_{\rho^*}^{\rho}\rangle d\rho^*(\theta).$$ Therefore, ${\bf{0}}$ is a subdifferential by definition of geodesic convexity. Thus, we have
\begin{equation}
    F(\rho)-F(\rho^*)\geq \int_{\Theta}\langle {\bf{0}}, T_{\rho^*}^{\rho}(\theta)-\theta\rangle d\rho^*(\theta)+\dfrac{\mu}{2}\mathcal{W}_2^2(\rho^*,\rho)=\dfrac{\mu}{2}\mathcal{W}_2^2(\rho^*,\rho)
\end{equation}
by definition of $\mu$-convexity along geodesics. Therefore, $$\mathcal{W}_2(\rho^*,\rho_n)\leq \sqrt{\dfrac{2}{\mu}F(\rho_0)-F(\rho^*)}\dfrac{1}{(1+\mu\xi)^n}.$$
\end{proof}

 \begin{proof}[{\bf Proof of Theorem \ref{Thm_inexact}}] 

Since $F$ is $(-L)$ geodesically convex, $F$ is $(-1/\xi)$ geodesically convex. Under Assumption \ref{assum_4},  since we assume  $\nabla\frac{\delta F}{\delta \rho}(\overline\rho_{n+1})\in \partial{F}({\overline\rho_{n+1}})$,
\begin{align*}
    F(\overline\rho_{n+1})-F(\overline\rho_n)&\leq -\xi\int_{\Theta} \left\langle\nabla \dfrac{\delta F}{\delta \rho}(\overline\rho_{n+1}),\dfrac{T_{\overline\rho_{n+1}}^{\overline\rho_n}-\mathbf{id}}{\xi} \right\rangle d\overline\rho_{n+1} + \dfrac{1}{2\xi}\mathcal{W}_2^2(\overline\rho_{n+1},\overline\rho_n)& \\
    &=-\xi \int_{\Theta} \left\langle \nabla \dfrac{\delta F}{\delta \rho}(\overline\rho_{n+1}),\nabla\dfrac{\delta F}{ \delta \rho}(\overline\rho_{n+1})+\beta_{n+1} \right\rangle d\overline\rho_{n+1} &\\&\qquad \qquad +\dfrac{1}{2\xi}\int_{\Theta} \left\Vert\xi\nabla\dfrac{\delta F}{ \delta \rho}(\overline\rho_{n+1})+\xi\beta_{n+1} \right\Vert^2d\overline\rho_{n+1} &\text{(by   (\ref{Eqn_inexact}))}\\
    &= -\dfrac{\xi}{2}{\int_{\Theta} \left\Vert \nabla \dfrac{\delta F}{\delta \rho}(\overline\rho_{n+1})\right\Vert^2d\overline\rho_{n+1}}+\dfrac{\xi}{2}\int_{\Theta} \left\Vert\beta_{n+1}\right\Vert^2d\overline\rho_{n+1} &\\
    &\leq -{\mu}{\xi} (F(\overline\rho_{n+1})-F(\rho^*))+\dfrac{\xi}{2}\epsilon_{n+1}. & \text{(by PL (\ref{Eqn_wass_PL}))}
\end{align*}
Thus, we have the one-step contraction with an error term,
\begin{equation}
\label{eq:one_step_error}
    (1+{\mu}{\xi})(F(\overline\rho_{n+1})-F^*)\leq (F(\overline\rho_{n})-F^*)+\dfrac{\xi}{2}\epsilon_{n+1}.
\end{equation}
% By [Lemma 14, \citep{yao2024wasserstein}], we set $A=\frac{1}{1+\mu/L}$, $B=\frac{L/2}{1+\mu/(L)}$.

Assume $\epsilon_{n}\leq \epsilon$ for every $n$,
\begin{align*}
    F(\overline\rho_n)-F(\rho^*)&\leq \dfrac{1}{(1+\mu\xi)^n}(F(\rho_0)-F(\rho^*))+\dfrac{1}{1+\mu\xi}\dfrac{1-\frac{1}{(1+\mu\xi)^n}}{1-\frac{1}{1+\mu\xi}}\dfrac{\xi\epsilon}{2}\\
     &\leq\dfrac{1}{(1+\mu\xi)^n}(F(\rho_0)-F(\rho^*))+\dfrac{\epsilon}{2\mu}.
\end{align*}
\end{proof}

\begin{proof}[{\bf Proof of Corollary \ref{Cor_inexact}}]
Set $A=\frac{1}{1+\mu\xi}$ and $B=\frac{\xi/2}{1+\mu\xi}$. By Equation (\ref{eq:one_step_error}) and \citep[Lemma 14]{yao2024wasserstein}, we have the following two cases.

(a) If $\epsilon_n\leq C_{exp}{\gamma}^{n}$ with $\gamma\in(0,1)$,
\begin{equation}
     F(\overline\rho_n)-F^* \leq A^n(F(\rho_0)-F^*)+\dfrac{BC_{exp}}{\vert A-\gamma\vert} \max\{A,\gamma\}^{n+1}.
\end{equation}
when $\gamma\neq A$. If $\gamma=A$, we define $\tilde \gamma=(1+\gamma)/2$, then we have $\epsilon_n\leq C_{exp}{\tilde\gamma}^{n}$. Thus,
\begin{equation}
     F(\overline\rho_n)-F^* \leq A^n(F(\rho_0)-F^*)+\dfrac{BC_{exp}}{\vert A-\tilde\gamma\vert} \max\{A,\tilde\gamma\}^{n+1}.
\end{equation}

(b) If $\epsilon_n\leq C_{poly}n^{-\zeta}$ with $\zeta>0$, there is a constant $C(\zeta,A)$,
\begin{equation}
    F(\overline\rho_n)-F^*\leq A^n(F(\rho_0)-F^*)+\dfrac{BC_{poly}C(\zeta,A)\zeta}{n^{\zeta}}.
\end{equation}
\end{proof}

\section{Further numerical experiments and discussions}
\label{app:more_exp_details}

\subsection{Particle methods and functional approximation}
\label{subsec:particle+functional_app}

 When applying particle methods to our entropy-regularized total objective of mean-field neural network (\ref{eqn:MFLD_obj}), the functional approximation objective (\ref{Eqn_prox_op_compute}) can be expressed as,\begin{equation}
\begin{aligned}
\label{Eqn_learn_otmap}
   T^{n+1} & = \arg\min_{T}\ \Big\{ \dfrac{1}{N}\sum_{i=1}^{N}l(\dfrac{1}{m}\sum_{j=1}^m \varphi(T(\theta_j^n),x_i),y_i)+\dfrac{\lambda}{m}\sum_{j=1}^{m}\left\Vert T(\theta_j^n)\right\Vert^2  \\
   & \qquad \qquad \qquad -\dfrac{\tau}{m}\sum_{j=1}^{m}\log|\det\nabla T(\theta_j^n)|+\dfrac{1}{2m\xi}\sum_{j=1}^{m}\left\Vert T(\theta_j^n)-\theta_j^n\right\Vert^2 \Big\},
\end{aligned}
\end{equation}
where the change of variable for entropy \citep{mokrov2021largescale} is utilized.
In our work, we specifically employ a shallow neural network to learn the optimal transport map $T^{n+1}$ at each iteration, using the right-hand side of (\ref{Eqn_learn_otmap}) as the loss function.

\subsection{Further discussions on MFLD experiments}
\label{app:further_exp}
  
\textbf{Further details of experiments.} The normalizing flow has two coupling blocks as defined in Remark \ref{ref:realnvp}. Both $s()$ and $t()$ are one-hidden-layer fully connected neural networks with 100 hidden dimensions and $softplus()$ activation. The two coupling layers are composed in an alternating pattern,
such that the components that are left unchanged in one coupling layer are updated in the next one. For the experiment presented in Figure \ref{Fig_inexact}, we extend the number of training iterations for learning the OT map to 300 to promote more stable convergence of the learned OT map.

\textbf{Particle discretization error.} In this experiment, we set $\tau=0.1$. Other settings are the same as Section \ref{subsec_mfld_exp}. As shown in Figure \ref{Fig_particle}, as $m$ increases, $\mathcal{W}_2(\rho_n^m,\tilde\rho_\tau^*)$ at convergence decreases for both algorithms. This experiment supports our conjecture that particle discretization error dominates the bias, if we assume $\tilde\rho_{\tau}^*$ is a good approximation of $\rho^*$. 

The particle discretization error is well studied for noisy gradient descent and a quantitative uniform-in-time propagation of chaos result for the MFLD has been established in \citep{chen2023uniformintimepropagationchaosmean}. Specifically, the "distance" between the finite particle dynamics and the mean-field dynamics converges at rate $O(1/m)$ for all $t>0$ under the uniform LSI condition. Note that even though the ``distance'' in theoretical analysis of uniform propagation of chaos is not simply defined to be the $\mathcal{W}_2$ distance between empirical measure of finite particle dynamics and absolutely continuous measure of mean-field dynamics, it is empirically observed that the $\mathcal{W}_2$ distance also demonstrates similar propagation-of-chaos property in Figure \ref{Fig_noisygd}. 
% The remaining challenge is how to theoretically characterize the particle discretization error of the proximal algorithm.

\begin{figure}[H]
\centering
    \begin{subfigure}[t]{0.48\textwidth}
        \includegraphics[width=\textwidth]{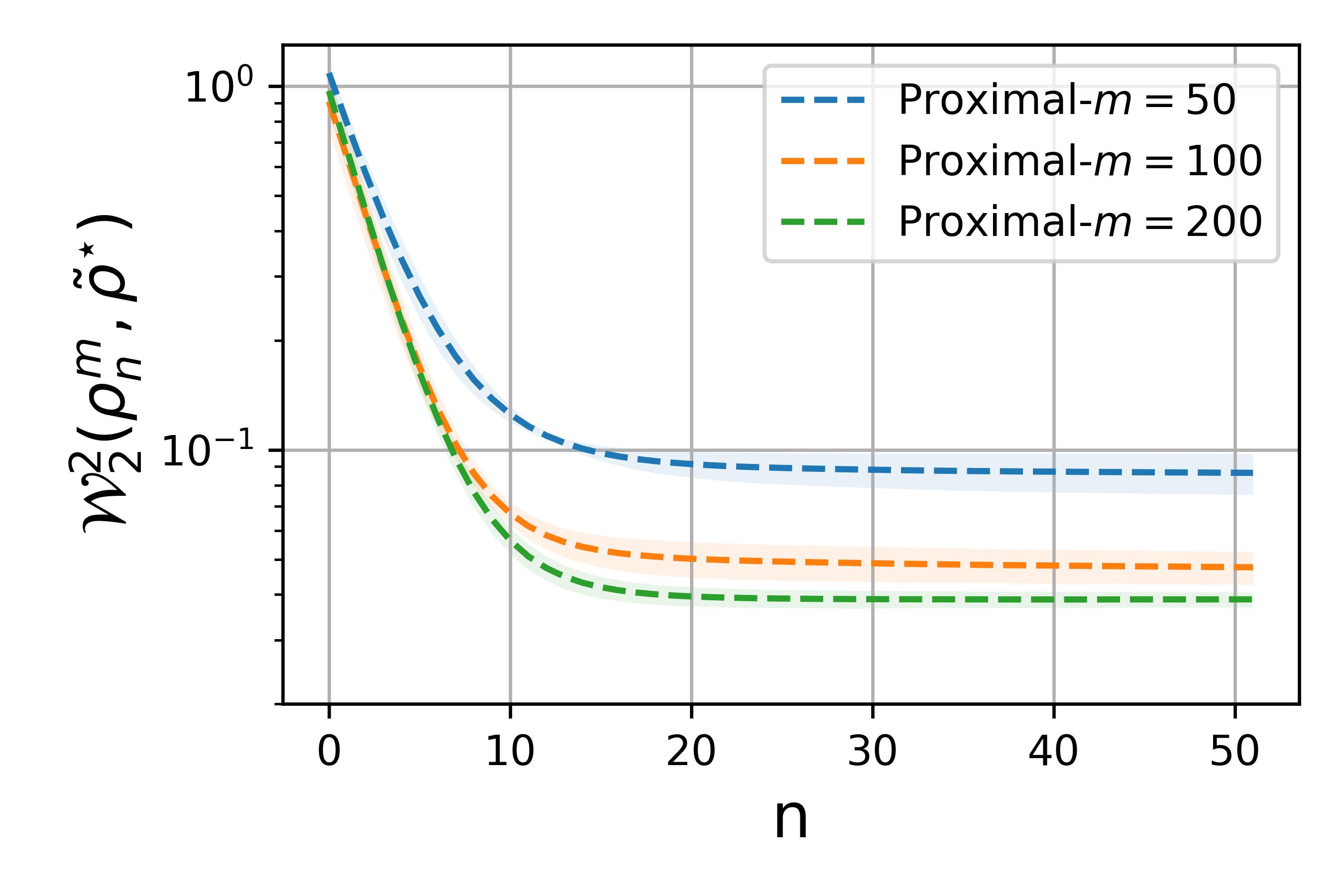}
        \caption{}
    \end{subfigure}
        \begin{subfigure}[t]{0.48\textwidth}
        \includegraphics[width=\textwidth]{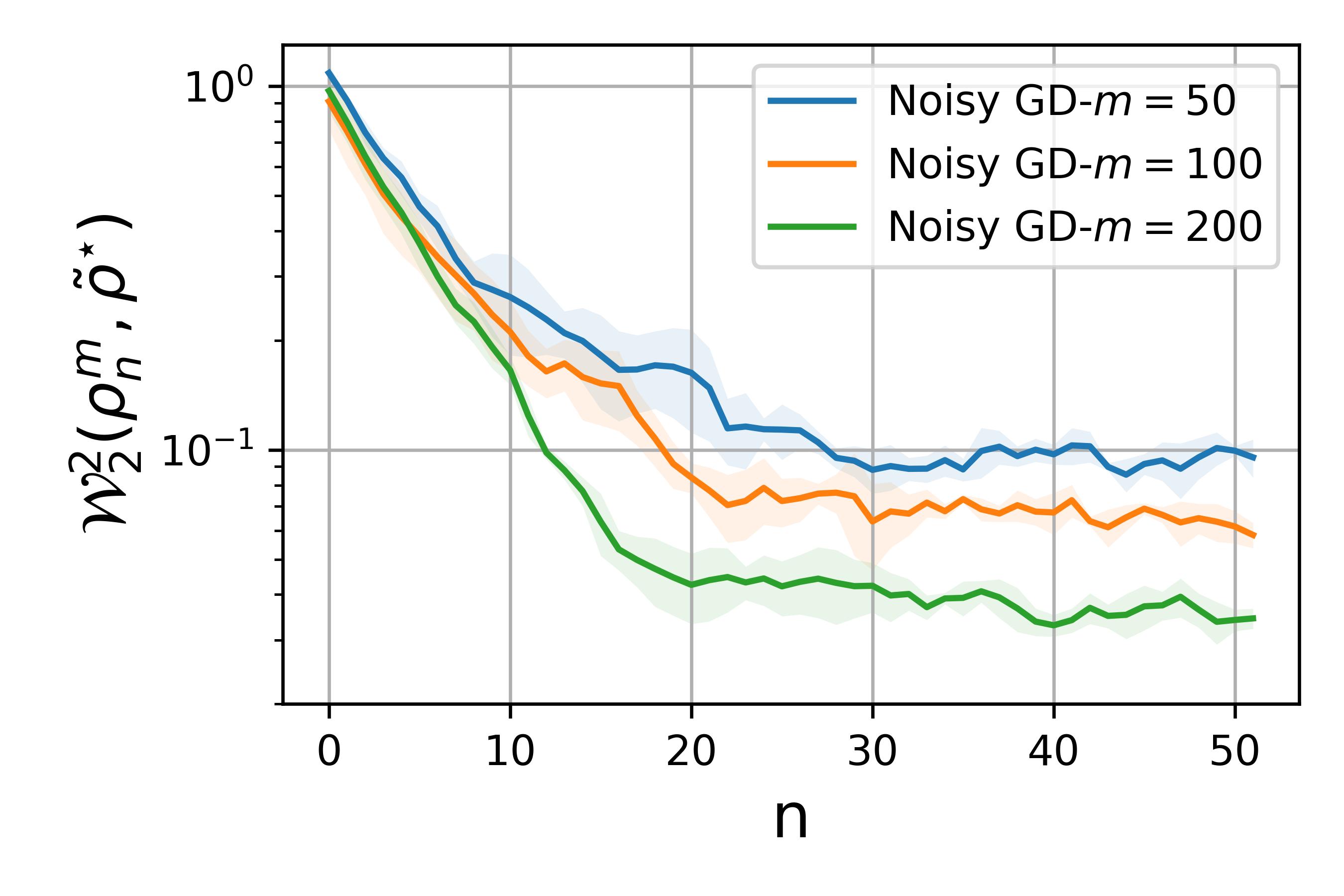}
        \caption{}
        \label{Fig_noisygd}
    \end{subfigure}
    \caption{Particle discretization error with different numbers of particles. We follow all the experiment settings in Section \ref{subsec_mfld_exp}.}
    \label{Fig_particle}
\end{figure}
\end{document}